\documentclass[sn-chicago]{sn-jnl}

\jyear{2022}%

\definecolor{wiasblue}   {cmyk}{1.0, 0.60, 0, 0}
\definecolor{mlugreen}{RGB}{172,6,52}

\def\E{\mathbb E}

\def\P{\mathbb P}

\def\R{\mathbb R}
\def\Z{\mathbb Z}
\def\mc{\mathcal}
\def\ms{\mathsf}

\def\lrsa{\leftrightsquigarrow}
\def\la{\lambda}

\def\su{\subseteq}

\def\de{\delta}

\def\Om{\Omega}

\def\d{{\rm d}}

\def\FF{\mc F}
\def\CC{\mc C}

\def\f{\frac}

\def\im{\item}
\def\sm{\setminus}
\def\th{\theta}
\def\bep{\begin{proof}}
\def\enp{\end{proof}}
\def\bepr{\begin{proposition}}
\def\enpr{\end{proposition}}
\def\bec{\begin{corollary}}
\def\enc{\end{corollary}}
\def\bea{\begin{align}}
\newcommand\eea{\end{align}}
\def\beas{\begin{align*}}
\def\eeas{\end{align*}}
\def\bet{\begin{theorem}}
\def\ent{\end{theorem}}
\def\bee{\begin{example}}
\def\ene{\end{example}}

\def\bede{\begin{definition}}
\def\ende{\end{definition}}
\def\ber{\begin{remark}}
\def\enr{\end{remark}}
\def\beca{\begin{cases}}
\def\enca{\end{cases}}
\def\bel{\begin{lemma}}
\def\enl{\end{lemma}}
\def\been{\begin{enumerate}}
\def\enen{\end{enumerate}}

\def\beit{\begin{itemize}}
\def\enit{\end{itemize}}
\def\befr{\begin{frame}}
\def\enfr{\end{frame}}

\renewcommand\le{\leqslant}
\renewcommand\ge{\geqslant}

\def\becb{\begin{tcolorbox}[colback=Dandelion!20]}
\def\encb{\end{tcolorbox}}

\def\bef{\begin{figure}[!h]}
\def\enf{\end{figure}}
\def\betp{\begin{tikzpicture}}
\def\entp{\end{tikzpicture}}
\def\co{\colon}
\def\endo{\end{document}}

\def\Inf{\ms{Inf}}
\def\Piv{\ms{Piv}}

\def\pa{\partial}

\renewcommand\it{\item}

\theoremstyle{plain}
\newtheorem{theorem}{Theorem}
\newtheorem{proposition}[theorem]{Proposition}
\newtheorem{corollary}[theorem]{Corollary}
\newtheorem{lemma}[theorem]{Lemma}

\theoremstyle{definition}
\newtheorem{definition}[theorem]{Definition}

\newtheorem{example}[theorem]{Example}

\theoremstyle{remark}
\newtheorem{remark}[theorem]{Remark}

\usepackage{dsfont}

\begin{document}

\title{Face and cycle percolation}

\author[1,2]{\fnm{Christian} \sur{Hirsch}}\email{hirsch@math.au.dk}

\author[3]{\fnm{Daniel} \sur{Valesin}}\email{daniel.valesin@warwick.ac.uk}

\affil[1]{\orgdiv{Department of Mathematics}, \orgname{Aarhus University}, \orgaddress{ \city{Aarhus},   \country{Denmark}}}
\affil[2]{\orgdiv{DIGIT Center}, \orgname{Aarhus University}, \orgaddress{ \city{Aarhus},   \country{Denmark}}}

\affil[3]{\orgdiv{Department of Statistics}, \orgname{University of Warwick}, \orgaddress{ \city{Coventry},   \country{United Kingdom}}}



\abstract{
	We consider face and cycle percolation as models for continuum percolation based on random simplicial complexes in Euclidean space. Face percolation is defined through infinite sequences of $d$-simplices sharing a $(d-1)$-dimensional face. In contrast, cycle percolation demands the existence of infinite $d$-cycles, thereby generalizing the lattice notion of plaquette percolation. We discuss the sharp phase transition for face percolation and derive comparison results between the critical intensities for face and cycle percolation. Finally, we consider an alternate version of simplex percolation, by declaring simplices to be neighbors whenever they are sufficiently close to each other, and prove a strict inequality involving the critical intensity of this alternate version and that of face percolation.
}

\maketitle
\section{Introduction}
\label{int_sec}

Rooted in statistical physics, percolation theory has developed into a vigorous research area with profound mathematical results and a variety of connections to different branches in the natural sciences. Although percolation theory has grown enormously in the last decades, there is still a reasonably clear consensus about its scope. Loosely speaking, it concerns existence, uniqueness and more refined properties of giant connected components in random graphs.

While random networks are designed to encode pairwise interactions, this modeling paradigm may be insufficient in order to capture complex phenomena characterized by higher-order interactions. In order to move beyond network models, simplicial complexes have become highly popular. This leads to the question of how to extend percolation theory into the domain of simplicial complexes. Henceforth, we describe two such candidates assuming that the reader has some familiarity with simplicial complexes in Euclidean space. In Section \ref{mod_sec}, we will review briefly the basic notions of computational topology that will be needed to carry out the proofs of our main results.

A natural first idea is the concept of \emph{$d$-face percolation}, which corresponds to the usual notion of percolation on the graph with vertex set formed by the $d$-simplices and where $d$-simplices are considered adjacent if they share a common $(d - 1)$-dimensional face. This adjacency notion has appeared in a variety of incarnations in the literature. Isolated vertices in the aforementioned notion of adjacency were considered both in the \v Cech and the Vietoris-Rips filtration in \citep{isol}. There, it is also noted that the Vietoris-Rips variant allows for a natural interpretation in terms of percolation of cliques, studied earlier in the context of Erd\H{o}s-R\'enyi graphs in \citep{clique}. The notion of face percolation  also makes sense for cubical complexes and has received particular attention in the case $d = D - 1$, see \citep{kozma,mikami}. In particular, when investigating Bernoulli thinnings of the $(D - 1)$-faces in $\Z^D$, the term \emph{plaquette percolation} has gained popularity.

Although not considered in the present work, we want to point out that \citep{isol} also discusses the concept of \emph{up-adjacency} where two $d$-simplices are \textit{adjacent} if they are contained in a common $(d + 1)$-simplex. For instance, when working with the \v Cech complex, this form of percolation is equivalent to percolation of the set of $d$-covered points. When working on a Poisson point process in $\R^D$,  recently the sharpness of the phase transition with respect to up-adjacency could be achieved in \citep{lpy} thereby extending earlier results on the existence of a phase transition for more general point processes \citep{bartekPerc1}. The sharpness is understood in the sense that below the critical intensity, connection probabilities decay exponentially fast in the distance.   

Our first main result, Theorem \ref{sharp_thm}, establishes that the phase transition for face percolation in the Poisson Vietoris-Rips complex is sharp. It is observed in \citep[Section 11]{lpy} that this form of percolation should be amenable to the method of continuous-time decision trees developed in that work. Nevertheless, we decided to present a proof relying on a discretization argument and  the classical OSSS method as developed in \citep{osss} for two reasons. First, in the setting of face percolation, the discretization does not introduce any subtle artifacts that would need to be removed through later arguments.  Hence, we do not need the full power of the machinery of continuous-time decision trees. Moreover, in the proof we highlight that in monotone Poisson-based models there is a very short argument relating the concept of influence with pivot probabilities. We believe that this connection is also of independent interest.

Although the previous paragraphs illustrate that face percolation is an important concept, it is not tailored to the problems arising in topological data analysis. More precisely, here, a central motive is the systematic investigation of the structure of components, loops and higher-dimensional features associated with simplicial complexes. One major obstacle in deriving central limit theorems that could form the basis for goodness-of-fit tests are the correlations induced by structures percolating through macroscopic regions in space. When thinking about the simplest topological feature, namely connected components, this corresponds precisely to the classical notion of continuum percolation. However, in order to create, for instance, large cavities with triangles in 3D, it does not suffice that the model be able to produce long sequences of adjacent triangles. Indeed, the triangles forming the boundary of such a cavity need to align along a surface without holes. Recently, homological percolation was proposed as an exciting generalization of classical continuum percolation to the context of simplicial complexes \citep{homPerc}. There are also close connections to embedding problems in lattice percolation \citep{lipschitz,lima}. However, all of these generalizations were not designed to capture the correlations that need to be controlled when deriving limit results.

%
%
For these reasons, we propose the concept  of \emph{$d$-cycle percolation}, where we require more than just the existence of an unbounded connected family of simplices, and we also exclude the occurrence of cavities.   Our notion of cycle percolation extends to the continuum the percolation analysis in \citep{kozma} of $(D - 1)$-dimensional plaquettes. To that end, we build on the notion of \emph{$d$-cycles} \citep{yvinec}. Phrased in plain language, a $d$-cycle is a collection of $d$-simplices where each $(d - 1)$-face lies in an even number of simplices in the collection. For instance, a finite 1-cycle is a union of loops, and an infinite connected 1-cycle contains a bi-infinite path. We say that there is \emph{$d$-cycle percolation} if there exists a connected infinite $d$-cycle.

%
%
Our second main result, Theorem \ref{pt_thm}, sheds light on the relations between the critical intensities of face and cycle percolation. Any infinite $d$-cycle contains an infinite sequence of adjacent $d$-faces, so that $d$-cycle percolation implies $d$-face percolation. However, loosely speaking, we can extract an infinite connected $(d - 1)$-cycle from the boundary of a sequence of adjacent $d$-faces, which gives us also an inequality in the other direction. Finally, it is also intuitive that if the intensity of points is so high that all connected components in the vacant phase of continuum percolation are bounded, then the occupied phase contains an infinite connected $d$-cycle.

%
%
From a topological perspective, it is highly plausible that the inequalities between the aforementioned critical intensities are strict. That is, there is a regime for the intensity where we observe an infinite sequence of adjacent $d$-faces but not yet an infinite connected $d$-cycle. We could not prove this statement. Instead, in the third main result, Theorem \ref{enhan_thm}, we give a strict inequality involving the critical intensity of face percolation and the critical intensity of an alternate model, which we call \textit{$*$-percolation of simplices}. In this model, simplices are declared adjacent whenever they are sufficiently close to each other (even if they do not share a face).

The rest of the manuscript is organized as follows. In Section \ref{mod_sec}, we define rigorously the concepts of face percolation, cycle percolation, and $*$-percolation of simplices. We also state the main results and provide quick hints on the main techniques used in the proof. Then, in Sections \ref{sharp_sec}, \ref{pt_sec} and \ref{sec:enhan}, we prove Theorems \ref{sharp_thm}, \ref{pt_thm} and \ref{enhan_thm}.

\section{Model and main results}
\label{mod_sec}
Let~$D \in \mathbb{N}$,~$d \in \{0,\ldots, D\}$ and~$x_1,\ldots,x_{d+1}\in \R^D$ be distinct points. In case~$d \ge 2$, assume that these points are \textit{affinely independent} (that is, they do not all lie in any affine subspace of~$\mathbb{R}^D$ with dimension less than~$d$). The convex hull of~$\{x_1,\ldots,x_{d+1}\}$ is called the~\textit{$d$-simplex} spanned by~$x_1,\ldots,x_{d+1}$; these points are called the \textit{vertices} of the~$d$-simplex. In case~$d \ge 1$, the convex hull of any subset with~$d$ elements of~$\{x_1,\ldots,x_{d+1}\}$ is a~$(d-1)$-simplex, and we refer to it as one of the \textit{faces} of the~$d$-simplex. 

Given a discrete set~$\Lambda \subset \R^D$, define the collection 
\begin{equation} \label{eq_scrS}
	\mathscr{S}_d(\Lambda):= \left\{ \begin{array}{c}
		\text{ConvexHull}(\{x_1,\ldots,x_{d+1}\}):x_1,\ldots,x_{d+1} \in \Lambda,\\[.1cm]
		\text{ all distinct with } \vert x_i - x_j \vert \le 1 \text{ for all }i,j
	\end{array}\right\},
\end{equation}
where~$\vert \cdot \vert$ denotes the Euclidean norm.

Let $X \su \R^D$, $D \ge 2$ be a homogeneous Poisson point process with intensity~$\la > 0$ defined on some probability space $(\Om, \FF, \P_\la)$ (we occasionally omit the subscript~$\lambda$ from~$\P_\lambda$). With probability one, for any~$d \ge 2$, any~$d+1$ distinct points~$x_1,\ldots,x_{d+1} \in X$ are affinely independent. Hence, for any~$d \in \{0,\ldots, D\}$, the collection~$\mathscr{S}_d(X)$ is a collection of simplices.  It is the collection of~$d$-simplices of~$X$ obtained from the \textit{Vietoris-Rips simplicial complex}; see Chapter~2 of~\citep{yvinec}. (Actually, it is common to have an extra parameter~$r > 0$ and the condition~$\vert x_i - x_j \vert \le r$ instead of~$\vert x_i - x_j \vert \le 1$ in the definition of~$\mathscr{S}_d(\Lambda)$. However, by scaling properties of the Poisson point process, there is redundancy in having both parameters~$\lambda$ and~$r$, so we normalize~$r = 1$).  Note that we take~$D$ as fixed throughout and omit it from the notation.

We now define the concept of face percolation announced in Section \ref{int_sec}.  For~$1 \le d \le D$, we say two $d$-simplices  are \textit{adjacent} if they share a face. We say that there is \emph{face percolation of~$d$-simplices} in~$X$ if~$\mathscr{S}_d(X)$ contains an infinite path of adjacent $d$-simplices, and we let
$$\lambda^\text{face}_d := \inf\big\{\la > 0\co \P_\la(\text{there is face percolation of~$d$-simplices in $X$}) > 0 \big\}$$
denote the corresponding {critical intensity} parameter.  In particular, $d = 1$ concerns the existence of an infinite path in classical continuum percolation.

Our first result is that the phase transition for face percolation is sharp. To make this precise, we let $\theta_{d,r}=\theta_{d,r}(\lambda)$ denote the probability of the event that there exists a sequence~$\sigma_1,\ldots,\sigma_n$ of $d$-simplices in $\mathscr{S}_d(X \cup \{o\})$ such that~$o \in \sigma_1$,~$\sigma_n$ intersects the complement of the ball~$B_r(o)$, and~$\sigma_i$ and~$\sigma_{i+1}$ are adjacent for each~$i$. Also let~$\theta_{d,\infty}= \theta_{d,\infty}(\lambda) := \lim_{r \to \infty} \theta_{d,r}(\lambda)$, which equals the probability that there is an infinite sequence of adjacent~$d$-simplices in~$\mathscr{S}_d(X \cup \{o\})$ starting at a simplex with a vertex equal to~$o$.

\bet[Sharp phase transition]
\label{sharp_thm} Let~$D \ge 2$ and~$d \in \{1,\ldots, D\}$. Then,
\begin{itemize}
	\item[(a)] for any~$\lambda < \lambda^\mathrm{face}_d$, we have
		\begin{equation*}
			\limsup_{r \to \infty} \;\frac{1}{r} \log \theta_{d,r}(\lambda) < 0;
		\end{equation*}
	\item[(b)] we have that
		\begin{equation*}
			\liminf_{\lambda \searrow \lambda^\mathrm{face}_d}\; \frac{\theta_{d,\infty}(\lambda)}{\lambda - \lambda^\mathrm{face}_d} > 0.
		\end{equation*}
\end{itemize}
\ent

%
%
For cycle percolation, we require more than just the existence of an unbounded connected family. Loosely speaking, we also exclude the occurrence of cavities, thus extending to the continuum the percolation analysis in \citep{kozma} of $(D - 1)$-dimensional plaquettes. 

Let~$\mathcal{M}$ be a collection of~$d$-simplices in~$\mathbb{R}^D$. We say that~$\mathcal{M}$ is a \textit{$d$-cycle} if for any~$(d-1)$-simplex~$\sigma$ with vertices in~$\R^D$, the number of~$d$-simplices of~$\mathcal{M}$ that have~$\sigma$ as a face is finite and even (possibly zero). We say that there is \textit{$d$-cycle percolation} in~$X$ if~$\mathscr{S}_d(X)$ contains an infinite~$d$-cycle~$\mathcal{M}$ that is \textit{face connected}, that is, any two simplices of~$\mathcal{M}$ are connected by a path of adjacent simplices of~$\mathcal{M}$. 
The corresponding critical parameter is defined as
$$\lambda^\mathrm{cycle}_d := \inf\big\{\la > 0\co \P_\la(\text{there is $d$-cycle percolation in $X$}) > 0 \big\}.$$

In order to state the next theorem, let us introduce the critical intensity for continuum percolation,~$\lambda^\mathrm{cont}(r)$, where~$r > 0$, as the infimum of the values of~$\lambda$ for which~$\cup_{x \in X}B_r(x)$ has an infinite connected component with positive probability. It is readily seen that face percolation of $1$-simplices coincides with continuum percolation with~$r=1$, so
\begin{equation*} \lambda^\mathrm{face}_1 = \lambda^\mathrm{cont}(1).\end{equation*}
Also let~$\lambda^\mathrm{vac}(r)$ as the supremum of the values of~$\lambda$ for which~$\R^D \backslash (\cup_{x \in X} B_r(x))$ has an infinite component with positive probability. 
\begin{theorem} \label{pt_thm}
For any~$D \ge 2$, we have that
	\begin{align*}
		&\lambda^\mathrm{cont}(1) = \lambda^\mathrm{face}_1 = \lambda^{\mathrm{cycle}}_1 \le \lambda^\mathrm{face}_2 \le \lambda^\mathrm{cycle}_2  \le \cdots \le \lambda^\mathrm{face}_{D-1} \le \lambda^\mathrm{cycle}_{D-1} \le \lambda^\mathrm{face}_D .
	\end{align*}
	Moreover, if~$D \ge 3$, we have
	\[\lambda^\mathrm{cycle}_{D-1} \le \lambda^\mathrm{vac}(1/2).\]
\end{theorem}

We conjecture that $\lambda^\mathrm{face}_d < \lambda^\mathrm{cycle}_d$ for all $2 \le d < D$, meaning that infinite face-connected $d$-cycles begin to appear at a strictly higher intensity than infinite paths of adjacent $d$-simplices.  In particular, by Theorem \ref{pt_thm} we would obtain also a strict inequality of critical intensities between the dimensions.

Our third theorem gives a strict inequality between the critical intensity for face percolation and the critical intensity for yet another notion of percolation, which we now define. Fix a range parameter~$r_0 > 0$ and say that two~$d$-simplices of~$\mathscr{S}_d(X)$ are \textit{$*$-adjacent} (with range~$r_0$) in case the Euclidean distance between them (seen as subsets of~$\R^D)$ is at most~$r_0$. We say that there is~\textit{$*$-percolation with range~$r_0$ of~$d$-simplices} in~$X$ if there is an infinite sequence~$\sigma_1,\sigma_2,\ldots$ of distinct~$d$-simplices of~$\mathscr{S}_d(X)$ such that~$\sigma_i$ and~$\sigma_{i+1}$ are~$*$-adjacent (with range~$r_0$) for every~$i$.
Define
\begin{align*}&\lambda^*_{d,r_0} \hspace{-.1cm}:= \inf\{\lambda\hspace{-.07cm} >\hspace{-.07cm} 0: \P_\lambda(\text{there is $*$-percolation with range~$r_0$ of $d$-simplices})\hspace{-.07cm}>\hspace{-.07cm} 0\}.\end{align*}

	\begin{theorem}\label{enhan_thm}
		For any~$D \ge 2$,~$d \in \{1,\ldots, D\}$ and~$r_0 > 2$, we have that
		\[\lambda^*_{d,r_0} < \lambda^\mathrm{face}_d.\]
	\end{theorem}
	The proof of this result uses percolation enhancement techniques, implemented by means of a spatial exploration of the point process, as in the approach developed in~\citep{franceschetti2011strict}. However, the nature of simplex percolation introduces several technical difficulties in following this approach. For instance, if the point process $X$ is revealed inside a region $A$ but not on~$A^c$, then there may be simplices of $\mathscr{S}_d(X)$ that are not yet revealed, but that involve vertices that have already been seen. Moreover, there are two effective range parameters in our setting: one of them, equal to 1, is the threshold for the formation of a $d$-simplex, and the other, $r_0$, is the threshold for $*$-adjacency. These complications make the enhancement argument much more delicate. While we believe that the above theorem should be true for any~$r_0 > 0$, we could only complete the proof with the assumption~$r_0 > 2$.

\section{Sharpness of phase transition: proof of Theorem \ref{sharp_thm}}
\label{sharp_sec}
The overall strategy is to adapt the OSSS methodology from \citep{osss}, which has recently been applied very successfully to establish the sharp phase transition in a variety of models, such as continuum percolation with unbounded radii \citep{raoufi_sub}. 

Since face percolation has a bounded range of dependence, the arguments leading to the exponential below a modified percolation threshold simplify substantially in comparison to \citep{raoufi_sub}. On the other hand, the geometric complexities coming from simplicial complexes imply that applying the OSSS method in the context of simplicial percolation requires a novel argument in order to relate influences to pivot probabilities.

In the entire section, we fix $d \ge 1$. In particular, we write $\th_r(\la)$ instead of $\th_{d,r}(\la)$. Moreover, we write~$\theta(\lambda)$ instead of~$\theta_\infty(\lambda)$.

Following \citep{raoufi_sub}, one important ingredient in the proof is to note that $d$-face percolation implies standard continuum percolation (with range 1). Moreover, for $\la < \lambda^\mathrm{cont}(1)$, i.e., in the sub-critical regime of continuum percolation with range 1, the probabilities to percolate beyond a certain distance decay exponentially in the distance.

To prove Theorem \ref{sharp_thm}, we will proceed as in \citep{raoufi_sub} and rely on the machinery of randomized algorithms developed in \citep{osss}. More precisely, the main step will be to establish the following differential inequality, see \citep[Lemma 1.7]{raoufi_sub}.

%
%
\bel[Differential inequality]
\label{sup_crit_lem}
For any~$\la_1,\la_2 > 0$ with~$\la_1 < \la_2$ there exists~$c_{\ms{DI}} > 0$ such that for any~$\la \in (\la_1,\la_2)$ and any~$r\ge 1$ we have
$$\f{\d}{\d\la}\log\th_r(\la) \ge {c}_{\ms{DI}}\cdot  \f{r}{\int_0^r \th_s(\lambda)\mathrm{d}s}.$$
\enl

Following  \citep{raoufi_sub}, Lemma \ref{sup_crit_lem} is the key in the proof of Theorem \ref{sharp_thm}. For completeness, we give the full argument.
\bep[Proof of Theorem \ref{sharp_thm}]
Define
\[S_r(\la):= \int_0^r \th_s(\la)\;\mathrm{d}s,\quad r > 0,\; \la > 0\]
and
\[\tilde{\la}:= \sup\left\{\la > 0:\; \limsup_{r\to \infty} \frac{\log S_r(\la)}{\log r} < 1\right\}.\]
Let us prove that~$\tilde{\lambda} \in (0,\infty)$. It is readily seen that~$\tilde{\la} \le \lambda^\mathrm{face}_d$. We claim that we also have~$\tilde{\la} \ge \lambda^\mathrm{cont}(1)$ (recall that~$\lambda^\mathrm{cont}(1)$ is the critical value for continuum percolation with range 1, which coincides with the critical value for percolation of 1-simplices).  To see this, first note that, if~$\mathcal{M}$ is a face connected set of~$d$-simplices of~$\mathscr{S}_d(X)$, the sets of vertices of the simplices in~$\mathcal{M}$ is contained in a cluster of continuum percolation. Then, use the well-known fact that in continuum percolation below criticality, the probability to percolate beyond a certain distance decays exponentially with the distance. This proves the claim.

Now fix~$\la_1,\la_2 > 0$ with~$\la_1 < \tilde{\la} <  \la_2$ and let~$c_{\ms{DI}}$ be a constant as in Lemma~\ref{sup_crit_lem}.  Take~$\lambda$,~$\lambda'$,~$\lambda''$ with~$\lambda_1 < \la < \la' < \la'' < \tilde{\la}$. Integrating the inequality in Lemma~\ref{sup_crit_lem} gives, for any~$r \ge 1$,
\begin{align*}
&\th_r(\la') \le \th_r(\la'') \cdot \exp\left\{-{c}_{\ms{DI}}\cdot (\la'' - \la')\cdot   \frac{r}{S_r(\la'')}\right\},\\
&\th_r(\la) \le \th_r(\la') \cdot \exp\left\{-{c}_{\ms{DI}}\cdot (\la' - \la)\cdot   \frac{r}{S_r(\la')}\right\}.
\end{align*}
Combining the first of these inequalities with the definition of~$\tilde{\la}$ and the fact that~$\la'' < \tilde{\la}$ shows that~$r \mapsto \th_r(\la')$ decays at least as fast as a stretched exponential function of~$r$. In particular,~$r \mapsto S_r(\la')$ is bounded. Using this with the second inequality above shows that~$r \mapsto \th_r(\la)$ decays exponentially.

To treat~$\la > \tilde{\la}$, we first give another definition and prove some auxiliary facts. Let
\[T_n(\la):= \frac{1}{\log n}\sum_{k=1}^n \frac{\th_k(\la)}{k},\quad n \ge 1,\; \la > 0,\]
and note that~$\lim_{n \to \infty} T_n(\la) = \th(\la)$. Moreover, for every~$\la \in (\la_1,\la_2)$ we bound:
\begin{align*}
	&\log n\cdot \frac{\d }{\d \la}T_n(\la) = \sum_{k=1}^n \frac{\frac{\d}{\d \la}\th_k(\la)}{k}  \ge {c}_{\ms{DI}}\sum_{k=1}^n \frac{\th_k(\la)}{S_k(\la)} \ge {c}_{\ms{DI}}\sum_{k=1}^n \frac{\int_k^{k+1} \th_s(\la)\; \d s}{S_k(\la)}\\
&\; = {c}_{\ms{DI}}\sum_{k=1}^n \frac{S_{k+1}(\la) - S_k(\la)}{S_k(\la)}  \ge {c}_{\ms{DI}}\sum_{k=1}^n \int_{S_k(\la)}^{S_{k+1}(\la)} \frac{\d t }{t} = {c}_{\ms{DI}}\cdot (\log S_{n+1}(\la) - \log S_1(\la)).\end{align*}
Hence,
\[ \frac{\d }{\d \la} T_n(\la) \ge \frac{{c}_{\ms{DI}}}{\log n} (\log S_{n+1}(\la) - \log S_1(\la)).\]

Now, let~$\la \in (\tilde{\la},\la_2)$, and also take~$\la' \in (\tilde{\la},\la)$. Integrating the above over $[\la', \la]$ and using monotonicity gives
\[T_n(\la) - T_n(\la') \ge \frac{{c}_{\ms{DI}}}{\log n}\cdot (\la - \la') \cdot (\log S_{n+1}(\la') - \log S_1(\la)).\]
Then,
\begin{align*}
	\th(\la) \ge \th(\la) - \th(\la') &= \lim_{n \to \infty} (T_n(\la) -T_n(\la')) \\
	&\ge {c}_{\ms{DI}}  (\la - \la')   \limsup_{n \to \infty} \frac{\log S_{n+1}(\la')}{\log n} \ge {c}_{\ms{DI}} (\la - \la') \ge {c}_{\ms{DI}} (\la - \tilde{\la}).
\end{align*}

It follows from what we have proved so far that~$\th(\la) = 0$ if~$\la < \tilde{\la}$ and~$\th(\la) > 0$ if~$\la > \tilde{\la}$; hence,~$\tilde{\la} = \lambda^\mathrm{face}_d$. This completes the proof.
\enp

%
%
In the rest of this section, we explain how to establish Lemma \ref{sup_crit_lem} via the OSSS technique. Let us introduce some notation. For shortness, we write $\th_r$ instead of $\th_r(\la)$. For~$x \in \R^D$ and~$r > 0$, we let~$Q_r(x):= x +[-r,r]^D$. For~$y \in X$ and~$A \subset \R^D$, we write~$y \lrsa A$ to denote the event that there exist~$d$-simplices~$\sigma_0,\ldots,\sigma_n$ in~$\mathscr{S}_d(X)$ such that~$y$ is a vertex of~$\sigma_0$,~$\sigma_n$ intersects~$A$, and~$\sigma_i,\sigma_{i+1}$ are adjacent for each~$i$. Note that~$\theta_r = \P(o \lrsa \partial B_r(o))$.

We now present an exploration algorithm that progressively reveals the clusters of adjacent~$d$-simplices intersecting the Euclidean sphere~$\pa B_s(o)$, with $s \le r$. The algorithm decides whether or not the event~$o \lrsa \partial B_r(o)$ occurs. We fix~$r > 0$,~$L \ge 2r$ and define the index set
$$I_L := \{x \in \Z^D \co \lvert x\rvert \le L\}.$$

\been
\im In a first phase, we reveal the Poisson process $X$ inside all the boxes $Q_1(x)$, $x \in \Z^D$ that intersect the set~$\{y \in \R^D: s-1 \le \lvert y\rvert \le s+1\}$.
\im Next, suppose that $X \cap Q_1(x_0), X \cap Q_1(x_1), \dots, X \cap Q_1(x_{t - 1})$ have already been revealed. Let $\CC_{t - 1}$ denote the union of all face connected components of revealed adjacent $d$-simplices that intersect the sphere $\partial B_s(o)$. Then, we reveal a box $Q_1(x_t)$ where $x_t \in I_L \sm\{x_0, \dots, x_{t - 1}\}$ is chosen such that $Q_1(x_t)$ is at distance at most 1 to an element of~$X$ that is a vertex of a simplex in $\CC_{t - 1}$.
\im If no such $x_t$ exists, the algorithm stops.
\enen

%
%
Then, the OSSS inequality reads as
\begin{align}
	\label{osss_eq}
\th_r (1 - \th_r) \le \sum_{x \in I_L} \de_s(x)\Inf(x),
\end{align}
where $\de_s(x)$ and $\Inf(x)$ are the \emph{revealment probability} and the \emph{influence} at position $x$, respectively.  That is, $\de_s(x)$ is the probability that in the exploration algorithm, the box $Q_1(x)$ is revealed; $\Inf(x)$ is the probability that when replacing the Poisson point process $X$ inside $Q_1(x)$ with an independent copy, then this changes whether or not the event $o \lrsa \partial B_r(o)$ occurs.

The key steps in the proof of Lemma \ref{sup_crit_lem} are now to give an upper bound on the revealment probability and an upper bound on the influence, the latter involving the derivative of the percolation probability. 
\bel[Revealment bound]
\label{reveal_lem}
For any~$x \in I_L$ we have
\[\int_0^r \de_s(x)\;\mathrm{d}s \le 4\sqrt{D} +  2\cdot 4^{D} \int_0^r \th_s\;\mathrm{d}s.\]
\enl
\bep
	Fix~$s \le r$ and~$x \in I_L$. In case~$\lvert \lvert x \rvert - s\rvert \le 2\sqrt{D}$, we simply bound~$\de_s(x) \le 1$. If on the other hand~$\lvert \lvert x \rvert - s\rvert >2\sqrt{D}$, then~$Q_1(x)$ is not revealed by the first step of the algorithm. Moreover, if~$Q_1(x)$ is revealed, then~$Q_2(x) = x + [-2,2]^D$ contains a node~$y \in X$ such that~$y \lrsa \partial B_s(o)$. Then, using the triangle inequality and Mecke's formula,
	\begin{align*}
		\de_s(x) &\le \P(\exists y \in X \cap Q_2(x): y \lrsa \partial B_s(o))\\[.18cm]
		&\le\P(\exists y \in X \cap Q_2(x): y \lrsa \partial B_{\lvert \lvert x \rvert -s \rvert - 2\sqrt{D}}(y))\\
		&\le \E[\#\{y \in X \cap Q_2(x): y \lrsa \partial B_{\lvert \lvert x \rvert -s \rvert - 2\sqrt{D}}(y)\}]= 4^D \cdot \th_{\lvert \lvert x \rvert -s \rvert - 2\sqrt{D}}.
	\end{align*}
The desired bound now follows from integrating over~$s$.
\enp

In order to give our influence bound, we will need to relate the influence at position $x$ to pivot probabilities. More precisely, we let $\Piv(x)$ denote the probability that adding a uniform point in $Q_1(x)$ changes whether or not the event $o \lrsa \partial B_r(o)$ occurs.

%
%
\bel[Influence bound]
\label{inf_lem}
For every $\la > 0$  and $x \in \Z^d$, we have
$\Inf(x) \le \la e^{\la}  \Piv(x)$.
\enl

\bep[Proof of Lemma \ref{inf_lem}]
First, the Poisson point process $X$ can be represented as
$$ X = (X\sm Q_1(x)) \cup\{Z_i\}_{i \le M},$$
where $M$ is a Poisson random variable with parameter $\la$ and $Z_1, Z_2, \dots$ are iid in $Q_1(x)$. Then, we let
$$N:=\inf_{n \ge 0}\big\{o \lrsa \partial B_r(o) \text{ in } X = (X\sm Q_1(x)) \cup\{Z_i\}_{i \le n}\big\}$$
denote the first index such that the event $\{o \lrsa \partial B_r(o)\} $ occurs in $X = (X\sm Q_1(x)) \cup\{Z_i\}_{i \le n}$, noting that $N$ may also take the values $0$ or $\infty$.

Next, we observe that the percolation model is increasing in the underlying Poisson point process $X$. In particular, if the event $\{o \lrsa \partial B_r(o)\}$ does not occur after resampling $X$ in $Q_1(x)$, then also $\{o \lrsa \partial B_r(o)\}$ does not occur in $X \sm Q_1(x)$. In particular, together with the independence of $M$ and $N$, we deduce that
$$
\Inf(x) \le \P(o \lrsa \partial B_r(o))  - \P(N = 0)\le  \P(M \ge N > 0) = \sum_{k \ge 1}\P(M\ge k) \P(N = k)
$$
Similarly, 
\[\Piv(x) = \P(N = M + 1) = \sum_{k \ge 1}\P(M = k - 1)\P(N = k).\]
The proof is concluded by noting that, for any~$k \ge 1$,
\begin{align*}
	\frac{\P(M \ge k)}{\P(M=k-1)} = \sum_{\ell = k}^\infty \frac{\lambda^\ell/\ell!}{\lambda^{k-1}/(k-1)!} = \sum_{\ell=k}^\infty \frac{\lambda^{\ell-k+1}}{\ell  (\ell-1)\cdots k} \le \la \sum_{j=0}^\infty \frac{\la^j}{j!}=\la e^{\la}.
\end{align*}
\enp

We now explain how to deduce Lemma \ref{sup_crit_lem} from Lemma~\ref{reveal_lem} and Lemma \ref{inf_lem}.

\bep[Proof of Lemma \ref{sup_crit_lem}]
Integrating over~$s$ in the OSSS inequality and using Lemma~\ref{reveal_lem} and Lemma \ref{inf_lem} gives
\[r \th_r (1 - \th_r) \le 2\la e^{\la}\left(4\sqrt{D} + 2\cdot 4^D\cdot  \int_0^r \th_s\;\mathrm{d}s\right)\sum_{x \in \Z^D} \Piv(x).\]
We now invoke Russo's formula to transform the sum of the pivot probabilities into the derivative and obtain
\[\f{\d}{\d\la}\th_r(\la) \ge \frac{r \th_r(\la)(1-\th_r(\la))}{ 2\la e^{\la}\left(4\sqrt{D} + 2\cdot 4^D\cdot\int_0^r \th_s\;\mathrm{d}s\right)}.\]

Now fix~$\la_1,\la_2 > 0$ with~$\la_1 < \la_2$. We claim that
\[ \inf\{1- \th_r(\la): \la \le \la_2,\; r \ge 1\} > 0.\] 
This follows from observing that if~$r \ge 1$, then~$1-\th_r(\la)$ is larger than or equal to the probability that the Euclidean ball of radius~$1$ centred at the origin contains no node of the Poisson point process. Next, we claim that
\[ \inf\left\{ \int_0^r \th_s\;\mathrm{d}s: \la \ge \la_1,\;r \ge 1\right\} > 0.\]
This is proved by noting that for~$r \in (0,1]$,~$\th_r(\la)$ is larger than the probability that there are~$d$ nodes in the Poisson point process inside~$B_1(o)$. Using these two bounds in the differential inequality obtained above allows us to obtain a constant~$c_{\ms{DI}}$ as required.
\enp

\section{Face and cycle threshold inequalities: proof of Theorem \ref{pt_thm}}
\label{pt_sec}

Let~$D \ge 1$ and~$d \in \{1,\ldots,D\}$. Let~$\mathcal{M}$ be a set of~$d$-simplices in~$\R^D$, and assume that~$\mathcal{M}$ is locally finite, in the sense that any ball of~$\R^D$ intersects finitely many simplices of~$\mathcal{M}$.  We define the \textit{simplex boundary of $\mathcal{M}$}, denoted~$\partial_\mathsf{s}\mathcal{M}$, as the set of~$(d-1)$-simplices
\[\partial_\mathsf{s}\mathcal{M}:= \{\pi:\;\#\{\sigma \in \mathcal{M}:\; \pi \text{ is a face of }\sigma\} \text{ is odd}\}.\]
We say that~$\mathcal{M}$ is a~${d}$-cycle in case~$\partial_\mathsf{s}\mathcal{M}$ is empty. For example, a 1-cycle is a set of 1-simplices (line segments) whose union is a collection of loops and (doubly infinite) lines. As another example: if a collection of two-dimensional triangles forming a surface in~$\R^3$ is a 2-cycle, then the surface has no cavities: the segments forming the boundary of a cavity would belong to the simplex boundary.

We can now give the first part of the proof of Theorem~\ref{pt_thm}.

\begin{proof}[Proof of Theorem~\ref{pt_thm}, $ \lambda^\mathrm{face}_1 = \lambda^\mathrm{cycle}_1$:]
	The first equality is obvious from the definitions of~$\lambda^\mathrm{cont}(1)$ and~$\lambda^\mathrm{face}_1$. The inequality~$\lambda^\mathrm{cycle}_1 \ge \lambda^\mathrm{face}_1$ is also obvious. For the reverse inequality we note that, when~$\lambda > \lambda^\mathrm{face}_1$, using well-known properties of continuum percolation, it is easy to see that there almost surely exists a doubly infinite sequence~$\ldots,\sigma_{-1},\sigma_0,\sigma_1,\ldots$ such that~$\vert \sigma_{i+1}- \sigma_i \vert \le 1$ for all~$i$. From such a sequence we obtain an infinite and face connected 1-cycle, completing the proof.
\end{proof}

\begin{proof}[Proof of Theorem~\ref{pt_thm}, $\lambda^\mathrm{face}_d \le \lambda^\mathrm{cycle}_{d+1}$:]
	Let~$d \ge 2$ and fix~$\lambda > \lambda^\mathrm{face}_d$. Inside an event of probability one, we can fix an infinite sequence of simplices~$\sigma_0,\sigma_1,\ldots \in \mathscr{S}_d(X)$ such that~$\sigma_{i}$ and~$\sigma_{i+1}$ share a face for each~$i$. Moreover, removing simplices from the sequence if necessary, we can assume that for each~$i$ and for each~$j > i+1$, there is no shared face between~$\sigma_i$ and~$\sigma_j$. For each~$i \ge 1$, we let~$\pi_i$ denote the face that is shared between~$\sigma_{i-1}$ and~$\sigma_i$.

	 Let~$\mathcal{M}$ denote the simplex boundary of~$\{\sigma_0,\sigma_1,\ldots\}$. Note that~$\mathcal{M}$ is equal to the set of faces of~$\sigma_0,\sigma_1,\ldots$ apart from~$\pi_1,\pi_2,\ldots$. In particular,~$\mathcal{M}$ is infinite and contained in~$\mathscr{S}_{d-1}(X)$. Moreover,~$\mathcal{M}$ is face connected, since any face of~$\sigma_{i-1}$  shares a ($(d-2)$-dimensional) face with any face of~$\sigma_i$. Lastly,~$\mathcal{M}$ is a~$(d-1)$-cycle. This follows from Proposition~11.5 in~\citep{yvinec}, which states that the simplex boundary of a simplex boundary is empty (to be precise, the proposition is stated in this reference for finite sets of simplices, but the proof, which just appeals to linearity module 2 of the simplex boundary with respect to the inclusion of simplices, applies equally well to the present context). 

	 Since~$\mathcal{M}$ is an infinite and face connected~$(D-1)$-cycle contained in~$\mathscr{S}_{D-1}(X)$, we obtain that~$\lambda > \lambda^\mathrm{cycle}_{D-1}$, so the proof is complete.
\end{proof}

To complete the proof of Theorem~\ref{pt_thm}, it remains to prove the inequality~$\lambda^\mathrm{cycle}_{D-1} \le \lambda^\mathrm{vac}(1/2)$. The strategy for doing so  will be to find, under the assumption that~$\lambda > \lambda^\mathrm{vac}(1/2)$, an infinite and face connected~$(D-1)$-cycle as the simplex boundary of a certain set of~$D$-simplices of the Delaunay tessellation~$\mathscr{D}$ obtained from the Poisson point process~$X$. We note that the Delaunay simplexes are the building blocks of the alpha-complex, which is a fundamental tool in computational topology \citep[Section 6.1]{yvinec}.

Recall that~$\mathscr{D}$ is the collection of~$D$-simplices defined by prescribing that the~$D$-simplex spanned by~$x_1,\ldots,x_{D+1} \in X$ belongs to~$\mathscr{D}$ if and only if the Voronoi cells of~$x_i$ and~$x_j$ have non-empty intersection for each~$i,j$. The collection~$\mathscr{D}$ is face connected, its simplices have disjoint interiors and their union is equal to~$\R^D$.

Since any face of a simplex of~$\mathscr{D}$ is shared with exactly one other simplex of~$\mathscr{D}$, the notion of simplex boundary of a subset of~$\mathscr{D}$ is particularly simple. In particular, given a set~$\mathcal{M} \subset \mathscr{D}$, denoting the topological boundary of a set by~$\partial$, we have
\[\partial \left(\bigcup_{\sigma \in \mathcal{M}}\sigma\right) = \bigcup_{\pi \in \partial_\mathsf{s}\mathcal{M}} \pi.\]
We define the \textit{simplex outer boundary} of~$\mathcal{M} \subset \mathscr{D}$, denoted~$\partial_\mathsf{s}^\mathrm{out}\mathcal{M}$, as the set of~$(D-1)$-simplices of~$\partial_\mathsf{s} \mathcal{M}$ that are shared faces between a simplex of~$\mathcal{M}$ and a simplex whose interior is in the unbounded connected component of~$\R^D \backslash (\cup_{\sigma \in \mathcal{M}}\sigma)$.

The following lemma about the simplex outer boundary of a set of simplices of~$\mathscr{D}$ will be useful. The proof is postponed to later in this section.

\begin{lemma}
	\label{lem_cyclePi}
	Let~$\mathcal{M}$ be a face connected and finite set of simplices of~$\mathscr{D}$. Then,~$\partial_\mathsf{s}^\mathrm{out} \mathcal{M}$ is a face connected~$(D-1)$-cycle.
\end{lemma}

Almost surely, each~$\sigma \in \mathscr{D}$ has a \textit{circumcenter}, that is a point~$x \in \R^D$ that is the center of the unique ball that has the vertices of~$\sigma$ in its boundary (and no other point of~$X$ in its closure). We write~$f(\sigma)$ to denote the circumcenter of~$\sigma$.

Let us say that~$\sigma \in \mathscr{D}$ is \textit{vacant} if~$f(\sigma) \notin \cup_{x \in X} B_{1/2}(x)$. The following observation, which follows readily from the triangle inequality, will be essential, so we record it as a remark.\\[-.2cm]

\begin{remark}\label{rem_vac}
	If~$\sigma\in \mathscr{D}$ is not vacant, then~$\sigma$ belongs to~$\mathscr{S}_D(X)$, and all the faces of~$\sigma$ belong to~$\mathscr{S}_{D-1}(X)$.
\end{remark}
$\;$\\[-0.6cm]

Let~$\sigma,\sigma' \in \mathscr{D}$ be adjacent (recall that this means that they share a face).  We say that~$\sigma,\sigma'$ are~\textit{$\mathsf{v}$-adjacent} if the line segment connecting~$f(\sigma)$ and~$f(\sigma')$ does not intersect~$\cup_{x \in X} B_{1/2}(x)$ (in particular,~$\sigma$ and~$\sigma'$ are both vacant).

We now define
\[\mathscr{K}_0:= \left\{ \begin{array}{c}
	\sigma \in \mathscr{D}:\; \sigma \text{ intersects the hyperplane}\\ \{(x_1,\ldots,x_D): x_1 = \cdots = x_{D-1} = 0\}
\end{array}
\right\}\]
and
\begin{equation}\label{eq_def_K}\mathscr{K}:= \mathscr{K}_0 \cup \left\{ \begin{array}{c}
	\sigma \in \mathscr{D}:\; \text{there exist } \sigma_1,\ldots, \sigma_n \in \mathscr{D}: \sigma_1 = \sigma, \sigma_n \in \mathscr{K}_0 \\
	\text{ and } \sigma_i,\sigma_{i+1} \text{ are $\mathsf{v}$ adjacent for each $i$}
\end{array}
\right\}.\end{equation}

\begin{lemma}
	\label{lem_unbounded}
	If~$\lambda > \lambda^\mathrm{vac}(1/2)$, then the set
	\[\mathbb{R}^D \backslash \left(\bigcup_{\sigma \in \mathscr{K}} \sigma\right)\] is almost surely unbounded.
\end{lemma}
We again postpone the proof of this lemma to later in this section, and now see how we can conclude the proof of Theorem~\ref{pt_thm}.

\begin{proof}[Proof of Theorem~\ref{pt_thm}, $\lambda^\mathrm{cycle}_{D-1} \le 2^D\lambda^\mathrm{vac}$.]
	Assume that~$\lambda > 2^D \lambda^\mathrm{vac}$. We will exhibit an infinite and face connected~$(D-1)$-cycle in~$\mathscr{S}_{D-1}(X)$, thus proving that~$\lambda > \lambda^\mathrm{cycle}(D-1)$. 

	Let~$\mathscr{K}$ be as in~\eqref{eq_def_K}.  Since~$\mathscr{K}$ is face connected, we can choose an enumeration~$\mathscr{K} = \{\sigma_1,\sigma_2,\ldots\}$ such that~$\mathscr{K}_n := \{\sigma_1,\ldots, \sigma_n\}$ is face connected for every~$n$.

By Lemma~\ref{lem_unbounded}, the topological outer boundary of~$\cup_{\sigma \in \mathscr{K}}\sigma$ is unbounded, hence the simplex outer boundary~$\partial^\mathrm{out}_\mathsf{s}\mathscr{K}$ is infinite. It is easy to see that for any~$\pi \in \partial_\mathsf{s}^\mathrm{out}\mathscr{K}$, we have~$\pi \in \partial_\mathsf{s}^\mathrm{out}\mathscr{K}_n$ for all~$n$ larger than some~$n_0$. Similarly, if~$\pi \notin \partial_\mathsf{s}^\mathrm{out}\mathscr{K}$, then~$\pi \notin \partial_\mathsf{s}^\mathrm{out}\mathscr{K}_n$ for all~$n$ larger than some~$n_0$. These facts are summarized by
	\begin{equation}\label{eq_limsets}\partial_\mathsf{s}^\mathrm{out} \mathscr{K} = \bigcup_{n\ge 1} \bigcap_{m \ge n} \partial_\mathsf{s}^\mathrm{out} \mathscr{K}_n = \bigcap_{n \ge 1}\bigcup_{m \ge n} \partial_\mathsf{s}^\mathrm{out} \mathscr{K}_n.\end{equation}
We can decompose~$\mathscr{F}:=\partial_\mathsf{s}^\mathrm{out} \mathscr{K}$ into face connected components; we let~$\mathscr{F}_0$ be one such component. By construction,~$\mathscr{F}_0$ is face connected. By Remark~\ref{rem_vac}, we have~$\mathscr{F}_0 \subset \mathscr{S}_{D-1}(X)$.

	Let us now prove that~$\mathscr{F}_0$ is a~$(D-1)$-cycle. Fix~$\pi \in \mathscr{F}_0$ and let~$\mu$ be a face of~$\pi$. By~\eqref{eq_limsets}, when~$n$ is large enough we have that~$\pi \in \partial_\mathsf{s}^\mathrm{out} \mathscr{K}_n$, and moreover, for any face~$\mu$ of~$\pi$,
	\[\{\pi' \in \partial_\mathsf{s}^\mathrm{out}\mathscr{K}_n: \; \mu \text{ is a face of }\pi'\}  = \{\pi' \in \partial_\mathsf{s}^\mathrm{out}\mathscr{K}: \; \mu \text{ is a face of }\pi'\}.\]
	By Lemma~\ref{lem_cyclePi},~$\partial_\mathsf{s}^\mathrm{out}\mathscr{K}_n$ is a~$(D-1)$-cycle, so the number on the left-hand side above is even. Hence,~$\mathscr{F}_0$ is a~$(D-1)$-cycle.

	 We will now prove that~$\mathscr{F}_0$ is infinite. Fix~$\bar{\pi} \in \mathscr{F}_0$. From the fact that~$\partial^\mathrm{out}_\mathsf{s}\mathscr{K}$ is infinite and~\eqref{eq_limsets}, it follows that the cardinality of~$\partial^\mathrm{out}_\mathsf{s} \mathscr{K}_n$ tends to infinity as~$n \to \infty$. Moreover, Lemma~\ref{lem_cyclePi} implies that~$\partial^\mathrm{out}_\mathsf{s} \mathscr{K}_n$ is face connected. Using these facts, it is possible to define increasingly longer chains of adjacent simplices of~$\partial^\mathrm{out}_\mathsf{s} \mathscr{K}_n$ starting at~$\bar{\pi}$. More precisely, for each~$n$ large enough that~$\bar{\pi} \in \partial^\mathrm{out}_\mathsf{s} \mathscr{K}_n$, we define a sequence~$\sigma_{n,0},\ldots,\sigma_{n,m_n}$ of distinct elements of~$\partial^\mathrm{out}_\mathsf{s} \mathscr{K}_n$ with~$\sigma_{n,0} = \bar{\sigma}$ and~$\sigma_{n,i}$ adjacent to~$\sigma_{n,i+1}$ for each~$i$, and so that~$m_n \to \infty$ as~$n \to \infty$. Combining this with~\eqref{eq_limsets} and a diagonal argument, we see that there is an infinite sequence~$\sigma_{\infty,0},\sigma_{\infty,1},\ldots$ of distinct elements of~$\partial^\mathrm{out}_\mathsf{s} \mathscr{K}_n$. This concludes the proof.
\end{proof}

We now turn to the two proofs of lemmas that we postponed, starting with Lemma~\ref{lem_cyclePi}.
\begin{proof}[Proof of Lemma~\ref{lem_cyclePi}.]
	Let us say that a set~$A \subset \R^D$ \textit{separates}~$\R^D$ if~$\R^D \backslash A$ has more than one connected component. It is sufficient to prove the lemma under the additional assumption that~$\cup_{\sigma \in \mathcal{M}} \sigma$ does not separate~$\R^D$. Indeed, if~$\cup_{\sigma \in \mathcal{M}}\sigma$ separates~$\R^D$, then we can let~$\mathcal{M}'$ denote the collection of~$\sigma \in \mathscr{D}$ whose interior lies in a bounded component of~$\R^D \backslash \cup_{\sigma \in \mathcal{M}}\sigma$, and note that~$\mathcal{M} \cup \mathcal{M}'$ is face connected, has the same simplex outer boundary as~$\mathcal{M}$, and~$\cup_{\sigma \in \mathcal{M} \cup \mathcal{M}'} \sigma$ does not separate~$\R^D$.

	With the additional assumption that~$\cup_{\sigma \in \mathcal{M}} \sigma$ does not separate~$\R^D$, we have~$\partial_\mathsf{s}^\mathrm{out}\mathcal{M} = \partial_\mathsf{s} \mathcal{M}$. As already noted, Proposition~11.5 in~\citep{yvinec} implies that the simplex boundary of any finite set of~$D$-simplices is a~$(D-1)$-cycle. It remains to prove that~$\partial_\mathsf{s} \mathcal{M}$ is face connected.

	We follow the approach in the proof of Theorem~7.3 in~\citep{grimmett2004random}. Assume for a contradiction that~$\partial_\mathsf{s} \mathcal{M}$ is not face connected, and let~$\mathcal{P}_1,\ldots,\mathcal{P}_n$ denote its face connected components. We claim that for each~$i$,~$\cup_{\pi \in \mathcal{P}_i} \pi$ does not separate~$\R^D$. To see this, fix~$\bar\pi \in \partial_\mathsf{s} \mathcal{M} \backslash \mathcal{P}_i$ and fix a point~$x\in \bar \pi \backslash (\cup_{\pi \in \mathcal{P}_i}\pi)$. Then, since~$\mathcal{M}$ is face connected, for any~$y$ in the (topological) interior of~$\cup_{\sigma \in \mathcal{M}} \sigma$, there is a topological path between~$y$ and~$x$ that lies entirely in the interior of~$\cup_{\sigma \in \mathcal{M}} \sigma$ except for its endpoint~$x$. Moreover, for any~$y$ in~$\R^D \backslash (\cup_{\sigma \in \mathcal{M}}\sigma)$, there is a topological path between~$y$ and~$x$ that does not intersect~$\cup_{\sigma \in \mathcal{M}}\sigma$ except for its endpoint~$x$. Using these considerations, the claim is proved.

	Now, the sets~$\cup_{\pi \in \mathcal{P}_i} \pi$ are closed, connected, do not separate~$\R^D$, and for~$i \neq j$, the intersection~$(\cup_{\pi \in \mathcal{P}_i} \pi) \cap (\cup_{\pi \in \mathcal{P}_j} \pi)$ is a union of simplices whose dimension is at most~$D-3$. By Theorem~11 of~\citep[paragraph 59, Section II]{kuratowski}, we have that~$\cup_i \cup_{\pi \in \mathcal{P}_i} \pi$ does not separate~$\R^D$, a contradiction.
\end{proof}

We now turn to the proof of Lemma~\ref{lem_unbounded}. Some preparatory lemmas will be required. Recall that for~$\sigma \in \mathscr{D}$,~$f(\sigma)$ denotes the circumcenter of the vertices of~$\sigma$. We denote by~$g(\sigma)$ the radius of the circumcenter, that is, the distance between~$f(\sigma)$ and any of the vertices of~$\sigma$.

\begin{lemma}
	\label{lem_good_g}
	For~$\ell$ large enough we have
	\begin{equation}
		\mathbb{P}\left(g(\sigma) < \sqrt{\ell} \text{ for all $\sigma \in \mathscr{D}$ intersecting $B_\ell(o)$}\right) > 1 - e^{-\ell^{1/4}}.
	\end{equation}
\end{lemma}
\begin{proof}
	We will show that the complement of the event in the probability on the left-hand side is contained in
	\[\{\exists  x \in B_{2\ell}(o):\;B_{\ell^{1/4}}(x) \cap X = \varnothing\}.\]
	Once this inclusion is proved, the desired inequality follows from routine concentration bounds for the Poisson point process.

	To prove the claimed inclusion, we assume that there exists~$\sigma \in \mathscr{D}$  that intersects~$B_\ell(o)$ and has~$g(\sigma) \ge \sqrt{\ell}$. The ball~$B_{g(\sigma)}(f(\sigma))$ has radius larger than~$\sqrt{\ell}$ and intersects~$B_{\ell}(o)$, since it contains~$\sigma$. Taking~$\ell$ large enough, there exists~$x \in B_{2\ell}$ such that~$B_{\ell^{1/4}}(x)$ is contained in the interior of~${B}_{g(\sigma)}(f(\sigma))$. By the definition of~$f(\sigma)$ and~$g(\sigma)$, we have that there is no point of~$X$ in the interior of~${B}_{g(\sigma)}(f(\sigma))$, so we obtain that~$B_{\ell^{1/4}}(x) \cap X = \varnothing$.
\end{proof}

\begin{lemma}\label{lem_far}
	Assume that~$\lambda > \lambda^\mathrm{vac}(1/2)$. Then, there exists~$\varepsilon > 0$ such that, for~$\ell$ large enough, we have~$\mathbb{P}(A) < \exp\{-\ell^{\varepsilon}\}$, where~$A$ is the event that there exist vacant simplices~$\sigma_1,\ldots, \sigma_n \in \mathscr{D}$ such that~$\sigma_1$ intersects~$B_\ell(o)$,~$\sigma_n$ intersects~$\mathbb{R}^D \backslash B_{5\ell}(o)$, and~$\sigma_i$ and~$\sigma_{i+1}$ are~$\mathsf{v}$-adjacent for each~$i$.
\end{lemma}
\begin{proof}
	We bound~$\mathbb{P}(A) \le \mathbb{P}(A \cap A') + \mathbb{P}((A')^c)$, where~$A'$ is the event that~$g(\sigma) < \sqrt{5\ell}$ for every~$\sigma \in \mathscr{D}$ that intersects~$B_{5\ell}(o)$. Lemma~\ref{lem_good_g} implies that~$\mathbb{P}((A')^c) < \exp\{-(5\ell)^{1/4}\}$. Next, note that on the event~$A \cap A'$, we have that~$f(\sigma_1) \in B_{2\ell}(o)$ and~$f(\sigma_n) \in \mathbb{R}^D \backslash B_{4\ell}(o)$. Moreover, by the definition of~$\mathsf{v}$-adjacency, the line segments~$[f(\sigma_i),f(\sigma_{i+1})]$ are contained in~$\R^D \backslash (\cup_{x \in X} B_{1/2}(x))$ for each~$i$. This shows that~$A \cap A'$ is contained in the event that~$\R^D \backslash (\cup_{x \in X} B_{1/2}(x))$ has a component connecting~$B_{2\ell}(o)$ and~$\mathbb{R}^D \backslash B^{4\ell}(o)$. By the exponential bounds on clusters of~$\R^D \backslash (\cup_{x \in X} B_{1/2}(x))$ in the subcritical regime developed in~\citep{raoufi_sub}, this has probability smaller than~$\exp\{-\ell^\varepsilon\}$ for some~$\varepsilon > 0$, completing the proof.
\end{proof}

\begin{proof}[Proof of Lemma~\ref{lem_unbounded}.] Assume that~$\lambda > \lambda^\mathrm{vac}(1/2)$. By Lemma~\ref{lem_far} and the Borel-Cantelli lemma, it is readily seen that there exists an almost surely finite random variable~$H$ such that the half line~$\{(0,\ldots,0,t): t \ge H\}$ is contained in the complement of~$\cup_{\sigma \in \mathscr{K}}\sigma$. This completes the proof.
\end{proof}

\section{Strict inequality between critical values: proof of Theorem~\ref{enhan_thm}}
\label{sec:enhan}

Throughout this section, we fix~$D \ge 2$,~$d \in \{1,\ldots,D\}$ and~$r_0 > 0$. For a discrete set~$\Lambda \subset \R^D$, recall the definition of~$\mathscr{S}_d(\Lambda)$ in~\eqref{eq_scrS}. We define
\[f_n^\mathrm{face}(\Lambda):= \mathds{1}\left\{ 
\begin{array}{c}
	\text{there exist }\sigma_0,\ldots,\sigma_k \in \mathscr{S}_d(\Lambda):\\ \sigma_0 \cap B_1(o) \neq \varnothing,\; \sigma_k \cap \partial B_n(o) \neq \varnothing,\\
	\sigma_i \text{ and } \sigma_{i+1} \text{ are adjacent for all }i
\end{array}
\right\}\]
and
\[f_n^*(\Lambda):= \mathds{1}\left\{ 
\begin{array}{c}
	\text{there exist }\sigma_0,\ldots,\sigma_k \in \mathscr{S}_d(\Lambda):\\ \sigma_0 \cap B_1(o) \neq \varnothing,\; \sigma_k \cap \partial B_n(o) \neq \varnothing,\\
	\mathrm{dist}(\sigma_i,\sigma_{i+1}) \le r_0 \text{ for all }i
\end{array}
\right\},\]
where~$\mathrm{dist}$ denotes Euclidean distance between sets, and~$\mathds{1}$ the indicator function.
	We then let
	\[\hat\theta^\mathrm{face}_n(\lambda) := \mathbb{E}_\lambda[f_n^\mathrm{face}(X)],\qquad \hat\theta^*_n(\lambda) := \mathbb{E}_\lambda[f_n^*(X)],\]
	where~$\mathbb{E}_\lambda$ is the expectation associated to a probability measure under which~$X$ is a Poisson point process with intensity~$\lambda$ on~$\R^D$. It is easy to see that
	\begin{equation}\label{eq_for_face}\lambda > \lambda^\mathrm{face}_d \quad \text{ if and only if }\quad  \lim_{n \to \infty} \hat \theta^\mathrm{face}_n(\lambda) > 0\end{equation}
	and 
		\begin{equation}\label{eq_for_star}\lambda > \lambda^*_{d,r_0} \quad \text{ if and only if }\quad  \lim_{n \to \infty} \hat \theta^*_n(\lambda) > 0.\end{equation}

			To prove Theorem \ref{enhan_thm}, our arguments follow the general strategy from \citep{strictIneq}. Given a set~$\Lambda \subset \R^D$, we say that~$x \in \Lambda$ is a \textit{special point} if there exist~$y, z \in \R^D$ such that
\begin{gather*}B_{1/8}(y) \subset B_1(x),\quad B_{1/8}(z) \subset B_1(x),\quad \mathrm{dist}(B_{1/8}(y),B_{1/8}(z))>1,\\[.2cm]
\#(\Lambda \cap B_{1/8}(y)) = \#(\Lambda \cap B_{1/8}(z)) = d,\\[.2cm] \Lambda \cap B_{r_0}(x) \backslash (\{x\} \cup B_{1/8}(y) \cup B_{1/8}(z)) = \varnothing.\end{gather*}
Let us motivate this definition. If~$x$ is a special point of~$\Lambda$, then the family of~$d$-simplices~$\mathscr{S}_d(\Lambda)$ has exactly two simplices that intersect~$B_{r_0}(x)$: one of them, denoted~$S_y$, has vertex set consisting of~$x$ together with~$d$ vertices in~$B_{1/8}(y)$, and the other, denoted~$S_z$, has vertex set consisting of~$x$ together with~$d$ vertices in~$B_{1/8}(z)$. Since~$\mathrm{dist}(B_{1/8}(y),B_{1/8}(z))>1$, these simplices have no face in common. These simplices do not share any vertex with any other simplex in~$\mathscr{S}_d(\Lambda)$ (which follows from~$S_y,S_z \subseteq B_1(x)$, the fact that there is no point of~$\Lambda$ inside~$B_{r_0}(x)$ apart from those that are contained in these simplices, and the assumption that~$r_0 > 2$). In particular, if~$x$ is deleted from~$\Lambda$, then, apart from the fact that~$S_y$ and~$S_z$ disappear, there is no effect for the~$d$-simplex face clusters of~$\Lambda$. Importantly, the deletion of~$x$ could affect a large~$*$-percolation cluster.  

We will use the fact that for any bounded set~$B \subset \R^D$, there exists a constant~$c_B > 0$ such that the number of special points of~$\Lambda$ inside~$B$ is at most~$c_B$ (uniformly in~$\Lambda$).

Given~$\Lambda \subset \R^D$ and~$\xi^1,\xi^2:\R^D \to \{0,1\}$, we define~$\Gamma(\Lambda,\xi^1,\xi^2): \Lambda \to \{0,1\}$ by letting, for each~$x \in \Lambda$,
\begin{equation}\label{eq_def_Gamma}[\Gamma(\Lambda,\xi^1,\xi^2)](x) = \begin{cases} \xi^1(x) \cdot \xi^2(x)& \text{if $x$ is a special point of $\Lambda$};\\  \xi^1(x)& \text{otherwise.}\end{cases}\end{equation}
Next, define the \textit{thinning} of~$(\Lambda,\xi^1,\xi^2)$ as the set
\begin{equation*}
	\mathcal{T}(\Lambda,\xi^1,\xi^2):= \{x \in \Lambda:\; [\Gamma(X,\xi^1,\xi^2)](x) = 1\}.
\end{equation*}
We let~$\mathbb{P}_{\lambda,p,q}$ be a probability measure under which a random triple~$(X,\Xi^1,\Xi^2)$ is defined, where~$X$ is a Poisson point process on~$\R^D$ with intensity~$1$ and~$\Xi^1,\Xi^2$ are random mappings from~$\R^D$ to~$\{0,1\}$, with~$\Xi^1 \sim \otimes_{x \in \R^D} \mathrm{Ber}(p)$ and~$\Xi^2 \sim \otimes_{x \in \R^D} \mathrm{Ber}(q)$, with~$X,\Xi^1,\Xi^2$ independent. We let~$\E_{\lambda,p,q}$ be the associated expectation operator. We will omit~$\lambda,p,q$ from the notation when they are clear from the context or unimportant. We now define
\begin{align*}&\Theta_n^*(\lambda,p,q):= \mathbb{E}_{\lambda,p,q}[f_n^*(\mathcal{T}(X,\Xi^1,\Xi^2))],\\[.2cm]& \Theta_n^\mathrm{face}(\lambda,p,q):= \mathbb{E}_{\lambda,p,q}[f_n^\mathrm{face}(\mathcal{T}(X,\Xi^1,\Xi^2))].\end{align*}

We are now ready to state:
\begin{proposition}\label{prop_diff_ineq}
	For any~$\lambda > 0$ and~$p,q \in (0, 1)$, there exists~$\varepsilon_0 = \varepsilon_0(\lambda,p,q)>0$ depending continuously on~$\lambda,p,q$ such that, for all~$n > 0$,
	\begin{equation*}
		\frac{\partial \Theta_n^*}{\partial q}(\lambda,p,q) > \varepsilon_0 \cdot \frac{\partial \Theta_n^*}{\partial p}(\lambda,p,q).
	\end{equation*}
\end{proposition}
We will use the following simple consequence of this proposition:
\begin{equation} \label{eq_half_quarter}
	\forall \lambda > 0,\; \forall p \in (0,1), \;\exists \delta > 0:\;  \Theta^*_n(\lambda,p-\delta,\tfrac{1}{2}) > \Theta^*_n(\lambda,p+\delta,\tfrac{1}{4}) \; \forall n > 0.
\end{equation}

\begin{proof}[Proof of Theorem~\ref{enhan_thm}]
	Fix~$\lambda > \lambda^*_{d,r_0}$. Choose~$\delta > 0$ corresponding to~$\lambda$ and~$p:= {\lambda_{d,r_0}^*}/{\lambda}$ in~\eqref{eq_half_quarter}. We then have~$(p-\delta)\lambda < \lambda^*_{d,r_0}$, and, by~\eqref{eq_for_star},
	\begin{equation} \label{eq_quick_aux} \lim_{n \to \infty} \hat\theta^*_n((p-\delta)\lambda) = 0.\end{equation}
	Moreover, for any~$n > 0$ we have
	\begin{align*}
		\hat\theta_n^\mathrm{face}((p+\delta)\lambda) = \Theta_n^\mathrm{face}(\lambda,p+\delta,0) &= \Theta_n^\mathrm{face}(\lambda,p+\delta,\tfrac14) \le \Theta^*_n(\lambda,p+\delta,\tfrac14)\\[.2cm]
	&\hspace{-1cm}\stackrel{\eqref{eq_half_quarter}}{<} \Theta^*_n(\lambda,p-\delta,\tfrac12)\le  \Theta^*_n(\lambda,p-\delta,1) = \hat\theta_n^*((p-\delta)\lambda).\end{align*}
	Together with~\eqref{eq_quick_aux}, this implies that~$\lim_{n \to \infty} \hat\theta^\mathrm{face}_n((p+\delta)\lambda) = 0$. On the other hand, if we had~$\lambda_c^\mathrm{face} = \lambda_c^*$, we would also have~$\lambda > \lambda_c^\mathrm{face}$, so~$\lim_{n \to \infty}\hat\theta_n^\mathrm{face}((p+\delta)\lambda) > 0$ by~\eqref{eq_for_face}, a contradiction.
\end{proof}

To prove Proposition~\ref{prop_diff_ineq}, we first relate the derivatives that appear there to probabilities of pivotality events. In order to define these, we introduce some more notation.
For any~$\xi:\R^D \to \{0,1\}$,~$x \in \R^D$ and~$i \in \{0,1\}$, define~$\psi_{x,i}(\xi):\R^D \to \{0,1\}$ as
\[[\psi_{x,i}(\xi)](y) = \begin{cases}
	i&\text{if } y = x;\\ \xi(y)&\text{otherwise.}
\end{cases}\]
Also define
\begin{align*}&\mathrm{Piv}^1_n(x) := \{f_n^*(\mathcal{T}(X\cup \{x\},\psi_{x,1}(\Xi^1),\Xi^2)) \neq f_n^*(\mathcal{T}({X}\cup \{x\},\psi_{x,0}(\Xi^1),\Xi^2))\},\\[.2cm]
	&\mathrm{Piv}^2_n(x) := \{f_n^*(\mathcal{T}(X\cup \{x\},\Xi^1,\psi_{x,1}(\Xi^2))) \neq f_n^*(\mathcal{T}({X}\cup \{x\},\Xi^1,\psi_{x,0}(\Xi^2)))\}.
\end{align*}
Note that for~$j=1,2$, the occurrence of~$\mathrm{Piv}^j_n(x)$ does not depend on~$\Xi^j(x)$.

We now state
\begin{lemma} \label{lem_first_pivot}
	For any~$\lambda > 0$ and~$p,q \in (0,1)$, there exists~$\varepsilon_0 = \varepsilon_0(\lambda, p, q)$ depending continuously on~$\lambda,p,q$ such that
	\[\mathbb{P}_{\lambda,p,q}(\mathrm{Piv}_n^{2}(x)) \ge \varepsilon_0 \cdot \mathbb{P}_{\lambda,p,q}(\mathrm{Piv}_n^{1}(x)).\]
\end{lemma}

Before we start the proof of this lemma, let us see how it allows us to conclude.
\begin{proof}[Proof of Proposition~\ref{prop_diff_ineq}]
As a straightforward application of Russo's formula, we have
	\[\frac{\partial \Theta_n^*}{\partial p}(\lambda,p,q) = \int_{\R^D}\P_{\lambda,p,q}(\mathrm{Piv}_n^1(x))\;\mathrm{d}x,\qquad \frac{\partial \Theta_n^*}{\partial q}(\lambda,p,q) = \int_{\R^D}\P_{\lambda,p,q}(\mathrm{Piv}_n^2(x))\;\mathrm{d}x.\]
	The statement of the proposition then readily follows from Lemma~\ref{lem_first_pivot}.
\end{proof}

One difficulty in proving Lemma~\ref{lem_first_pivot} is that, given a triple~$(\Lambda,\xi^1,\xi^2)$, determining whether or not a vertex~$y \in \Lambda$ belongs to~$\mathcal{T}(\Lambda,\xi^1,\xi^2)$ depends not only on~$\xi^1(y)$ and~$\xi^2(y)$, but also on the presence and positions of other points of~$\Lambda$ near~$y$, which determine whether or not~$y$ is special. We now work towards overcoming this obstacle by switching to different pivotality events, which will allow us to ignore the ``special'' status of points, at least inside a certain region.

For any~$\xi: \R^D \to \{0,1\}$,~$x \in \R^D$ and~$\ell > 0$, define~$\varphi_{x,\ell}(\xi): \R^D \to \{0,1\}$ as
\[[\varphi_{x,\ell}(\xi)](y) = \begin{cases}
	1&\text{if } y \in B_{\ell+r_0}(x) \backslash \{x\};\\ \xi(y)&\text{otherwise.}
\end{cases}\]
Then let
\begin{align*}&\mathrm{Piv}^1_{n,\ell}(x) \\&:= \{f_n^*(\mathcal{T}(X\cup \{x\},\psi_{x,1}(\Xi^1),\varphi_{x,\ell}(\Xi^2))) \neq f_n^*(\mathcal{T}({X}\cup \{x\},\psi_{x,0}(\Xi^1),\varphi_{x,\ell}(\Xi^2)))\},\\[.2cm]
	&\mathrm{Piv}^2_{n,\ell}(x) \\
	&:= \{f_n^*(\mathcal{T}(X\cup \{x\},\Xi^1,\psi_{x,1}(\varphi_{x,\ell}(\Xi^2)))) \neq f_n^*(\mathcal{T}({X}\cup \{x\},\Xi^1,\psi_{x,0}(\varphi_{x,\ell}(\Xi^2))))\}.
\end{align*}
In words, these are events that~$x$ is pivotal as defined earlier, except that before considering this pivotality, we replace the triple~$(X,\Xi^1,\Xi^2)$ by $(X,\Xi^1,\varphi_{x,\ell}(\Xi^2))$, that is, we change~$\Xi^2$ so that it becomes identically equal to 1 in~$B_{\ell+r_0}(x) \backslash \{x\}$.

We relate the probabilities of the new pivotality events to those of the original ones through the following.
\begin{lemma}\label{lem_relate_piv}
	For any~$\ell > 0$ and~$p,q \in (0,1)$, there exist~$\underline{\alpha} = \underline{\alpha}(\ell, p,q) > 0$ and~$\bar{\alpha} = \bar{\alpha}(\ell,p,q) > 0$ such that for any~$\lambda >0$,~$n > 0$ and~$x \in \R^D$ we have 
	\[\underline{\alpha} \cdot  \P_{\lambda,p,q}(\mathrm{Piv}^j_n(x)) \le \P_{\lambda,p,q}(\mathrm{Piv}^j_{n,\ell}(x))  \le\bar{\alpha}\cdot  \P_{\lambda,p,q}(\mathrm{Piv}^j_n(x)),\quad j \in \{1,2\}.\]
	Moreover,~$\underline{\alpha}$ and~$\bar{\alpha}$ can be taken to depend continuously on~$p,q$.
\end{lemma}
\begin{proof}
We give the proof for~$j = 1$ ($j=2$ is treated in the same way).
	Recall the definition of~$\Gamma(\Lambda,\xi^1,\xi^2)$ in~\eqref{eq_def_Gamma}.
	The occurrence of~$\mathrm{Piv}^1_n(x)$  only depends on~$X$ and on~$\Gamma(X \cup \{x\},\Xi^1,\Xi^2)$, that is, we can write, for an appropriate function~$g$,
	\[\mathds{1}_{\mathrm{Piv}^1_n(x)} = g(X,\Gamma(X\cup \{x\},\Xi^1,\Xi^2)).\]
	 We then also have
	\[\mathds{1}_{\mathrm{Piv}^1_{n,\ell}(x)} = g(X,\Gamma(X\cup \{x\},\Xi^1,\varphi_{n,\ell}(\Xi^2))).\]
	We thus have
	\[\mathbb{P}_{\lambda,p,q}(\mathrm{Piv}^1_n(x)) = \int g\;\mathrm{d}\mu,\qquad \mathbb{P}_{\lambda,p,q}(\mathrm{Piv}^1_{n,\ell}(x)) = \int g\;\mathrm{d}\nu,\]
	where~$\mu$ is the distribution of~$(X,\Gamma(X\cup \{x\},\Xi^1,\Xi^2))$ under~$\P_{\lambda,p,q}$ and~$\nu$ is the distribution of~$(X,\Gamma(X\cup \{x\},\Xi^1,\varphi_{n,\ell}(\Xi^2)))$ under~$\P_{\lambda,p,q}$. Now, letting~$N(\ell)$ denote the maximum number of special vertices that can appear inside~$B_{\ell+r_0}(x)$, it is easy to see that the Radon-Nikod\`ym derivative~$\frac{\mathrm{d}\mu}{\mathrm{d}\nu}$ is bounded from below by~$\min(p,q)^{N(\ell)}$ and bounded from above by~$\min(p,q)^{-N(\ell)}$. This completes the proof.
\end{proof}
Lemma~\ref{lem_first_pivot} follows from Lemma~\ref{lem_relate_piv} and the following:
\begin{lemma}
	For~$\ell$ large enough, any~$\lambda > 0$ and~$p,q \in (0,1)$, there exists~$\varepsilon_1 = \varepsilon_1(\ell,\lambda,p,q)$ depending continuously on~$\lambda$,~$p$,~$q$ such that for any~$n >2\ell$ and~$x \in B_n(o)$, we have
	\begin{equation*}
		\mathbb{P}(\mathrm{Piv}^{2}_{n,\ell}(x)) \ge \varepsilon_1 \cdot \mathbb{P}(\mathrm{Piv}^{1}_{n,\ell}(x)).
	\end{equation*}
\end{lemma}
The strategy of proof of this lemma will be as follows. We will define an exploration process of the configuration~$(X,\Xi^1,\varphi_{x,\ell}(\Xi^2))$, starting with the complement of the ball~$B_{\ell}(x)$, and advancing radially towards the center of this ball. The exploration will fail at any time that it becomes impossible, given the information already revealed, for either of the events~$\mathrm{Piv}^{1}_{n,\ell}(x)$ or~$\mathrm{Piv}^{2}_{n,\ell}(x)$ to occur (so that~$\{\text{Exploration fails}\} \subset (\mathrm{Piv}^1_{n,\ell}(x) \cup \mathrm{Piv}^2_{n,\ell}(x))^c$). If such a failure does not happen, there will be three ``success'' scenarios where the exploration will be terminated (still leaving an unexplored region). These three scenarios will be mutually exclusive, and unless the exploration fails, one of them will necessarily occur. In particular, letting~$A_1$,~$A_2$,~$A_3$ denote the events that the exploration terminates in scenario 1, 2 or 3, respectively, we will have
\begin{equation}
	\label{eq_mutual}
\mathrm{Piv}^{1}_{n,\ell}(x)\cup \mathrm{Piv}^{2}_{n,\ell}(x) \subset \{\text{Exploration suceeds}\} = A_1 \cup A_2 \cup A_3. \end{equation}
 Let~$\mathcal{G}$ denote the~$\sigma$-algebra generated by the exploration up to when it terminates (either in a failure or in a success scenario). We will prove that there  exists~$\varepsilon_1 > 0$ such that
\begin{equation}
	\label{eq_scenarios}
	 \mathbb{P}(\mathrm{Piv}^{2}_{n,\ell}(x)\mid \mathcal{G}) \cdot \mathds{1}_{A_j} \ge \varepsilon_1\cdot \mathds{1}_{A_j}, \quad j = 1,2,3.
\end{equation}

Once these statements are proved, we can conclude with
\begin{align*}
	\mathbb{P}(\mathrm{Piv}^{2}_{n,\ell}(x)) &\stackrel{\eqref{eq_mutual}}{=} \sum_{i=1}^3 \P(\mathrm{Piv}^2_{n,\ell}(x) \cap A_i) \\
	&= \sum_{i=1}^3  \mathbb{E}\left[ \mathbb{P}(\mathrm{Piv}^{2}_{n,\ell}(x)  \mid \mathcal{G})\cdot \mathds{1}_{A_i}\right]\\& \stackrel{\eqref{eq_scenarios}}{\ge} \varepsilon_1 \sum_{i=1}^3 \P(A_i) =  \varepsilon_1 \cdot \mathbb{P}(A_1 \cup A_2 \cup A_3) \stackrel{\eqref{eq_mutual}}{\ge} \varepsilon_1 \cdot \P(\mathrm{Piv}^{1}_{n,\ell}(x)).
\end{align*}

It is worth pointing out that formally, this type of argument involving a ``spatial Markov property'' for Poisson point processes relies on the notion of \textit{stopping sets}, as in~\citep{lpy}, but we will refrain from introducing the underlying measure-theoretic notions here, and will simply describe the spatial construction in words.

We will describe the rest of the proof in this section for the case~$D=2$, since the geometric explanations involved in the proof become more clumsy in higher dimension, but no complication is added.

We fix~$\ell > 0$ and~$\eta > 0$. In what follows, we will assume that~$\ell$ is sufficiently large and~$\eta$ is sufficiently small. Also fix~$n > 2 \ell$ and~$x \in B_n(o)$.  We assume that~$2\ell < \vert x\vert < n - 2\ell$. The cases where~$\vert x\vert \le 2\ell$ and~$n-2\ell\le \vert x\vert \le n$ must be treated separately, but since they are similar to~$2\ell < \vert x \vert < n-2\ell$, only simpler, we will omit the details.

We let
\[\mathcal{X} := \mathcal{T}(X \cup \{x\},\Xi^1,\varphi_{x,\ell}({\Xi}^2)).\]
It will be useful to note that
\begin{equation} \label{eq_xcal}
	\mathcal{X} \cap (B_{\ell}(x)\backslash \{x\}) =  \{y \in X: \Xi^1(y) = 1\} \cap (B_{\ell}(x)\backslash \{x\}).
\end{equation}
Also let
\[ \mathcal{Y}_t := \mathcal{X} \cap \overline{(B_{\ell-t}(x))^c},\quad t \in [0,\ell).\]
(the overline denotes topological closure).

Note that by~\eqref{eq_xcal}, for each~$t \in (0, \ell)$ we have~$\mathcal{Y}_t = \mathcal{Y}_{t-} \cup (\{y \in X: \Xi^1(y)=1\} \cap \partial B_{\ell-t}(x))$.\\[-.2cm]

\noindent \textbf{Definition of exploration process.} To initialize the process, we reveal~$\mathcal{Y}_0$.  Consider  $*$-percolation clusters induced from~$\mathcal{Y}_0$. We define~$A_0$ as the event that none of these clusters intersect both~$B_1(o)$ and~$(B_n(o))^c$, but there is a (necessarily unique) cluster intersecting both~$B_1(o)$ and~$B_{\ell+r_0}(x)$, and at least one cluster intersecting both~$(B_n(o))^c$ and~$B_{\ell+r_0}(x)$. If~$A_0$ does not occur, we already declare the exploration as a failure and terminate it.

On~$A_0$, define~$U_o$ as the set of vertices of~$\mathcal{Y}_0$ belonging to the simplices of the cluster connecting~$B_1(o)$ to~$B_{\ell+r_0}(x)$, and define~$U_n$ as the set of vertices of~$\mathcal{Y}_0$ belonging to the simplices of all of the clusters connecting~$(B_n(o))^c$ to~$B_{\ell+r_0}(x)$. 

We now reveal~$\mathcal{X}$ radially towards~$x$, stopping when a newly revealed vertex causes either of~$U_o$ or~$U_n$ to grow. Formally, this is done as follows. For each~$t \in [0,\ell)$,~$U_o$ is contained in the set of vertices of the simplices of one of the clusters induced by~$\mathcal{Y}_t$; let~$U_o(t)$ denote the set of vertices of the simplices of this cluster. Similarly, we let~$U_n(t)$ denote the set of vertices of~$\mathcal{Y}_t$ belonging to the simplices of all of the clusters induced by~$\mathcal{Y}_t$ that intersect~$U_n$. In particular,~$U_o(0) = U_o$ and~$U_n(0) = U_n$. Next, still on the event~$A_0$, define
\[\tau:= \inf\{t \in (0,\ell):\;U_o(t) \neq U_o \text{ or } U_n(t) \neq U_t\},\]
with the convention that~$\inf \varnothing = \infty$ (note that either~$\tau = \infty$ or~$\tau \le r_0+1$).
Define the event
\[A_0':= A_0 \cap \{\tau < \infty\} \cap \{U_o(\tau) \neq U_n(\tau)\}.\]
If~$A_0 \cap (A_0')^c$ occurs, we again declare failure.

On~$A_0'$, there is a newly found vertex~$Y \in \mathcal{X} \cap \partial B_{\ell - \tau}(x)$ that caused exactly one of the sets~$U_o(\tau-)$ or~$U_n(\tau-)$ to grow. It will be convenient to re-label the sets~$U_o(\tau)$ and~$U_n(\tau)$ so as to be able to refer more easily to whichever of the two contains~$Y$; with this in mind, we define
\[V_1:= \begin{cases} U_o(\tau)&\text{if } Y \in U_o(\tau),\\ U_n(\tau)&\text{if } Y \in U_n(\tau);\end{cases} \qquad V_2:= \begin{cases} U_n(\tau)&\text{if } Y \in U_o(\tau),\\ U_o(\tau)&\text{if } Y \in U_n(\tau).\end{cases}\]
It is worth noting that
	\begin{equation} \label{eq_v2_empty} V_2 \cap\overline{B_{\ell}(x) \backslash B_{\ell-\tau}(x)}= \varnothing. \end{equation}
We also define, on the event~$A_0'$, the regions
	\[\mathcal{R}_i := \bigcup_{y \in V_i} B_{r_0}(y),\quad i = 1,2.\]

We now define~$A_1$, the event corresponding to the first of the three successful termination scenarios for the exploration, as
\[A_1 := A_0' \cap  \{ V_2 \cap B_{\ell - \tau + \eta}(x) \neq \varnothing\}\]
(recalling that~$\eta$ will be taken sufficiently small). Note that, by~\eqref{eq_v2_empty},~$A_1$ can only occur if~$\tau <\eta$ (so that~$\ell - \tau + \eta > \ell$, and a point of intersection between~$V_2$ and~$B_{\ell - \tau + \eta}(x)$ could then exist in~$B_{\ell - \tau +\eta}(x)\backslash B_{\ell}(x)$). 

In the event~$A_0' \cap A_1^c$, we skip (for now) the annulus~$B_{\ell - \tau}(x) \backslash B_{\ell - \tau - \eta}(x)$, and we explore the region~$B_{\ell - \tau - \eta}(x) \cap \mathcal{R}_2 \cap \mathcal{R}_1^c$ radially towards~$x$, stopping either if we find a vertex of~$\mathcal{X}$ or we reach~$\partial B_{\ell - \tau - r_0 + \eta}(x)$ without finding a vertex of~$\mathcal{X}$. More formally, define
\[\sigma := \inf\{t \in (\eta,\;r_0 - \eta]: \mathcal{X} \cap \partial B_{\ell - \tau - t}(x) \cap \mathcal{R}_2 \cap \mathcal{R}_1^c  \neq \varnothing\},\]
again with the convention that~$\inf \varnothing = \infty$. In case~$\sigma < \infty$, let~$Z^*$ be the newly found vertex of~$\mathcal{X}$. Note that the appearance of~$Z^*$ could cause~$V_2$ to increase by a simplex (or multiple simplices). It could even cause~$V_1$ to increase, because a simplex might be formed between~$Z^*$ and vertices in~$\mathcal{R}_1$ near the boundary of~$\mathcal{R}_1$. If~$Z^*$ causes~$V_1$ to increase, then this causes the updated versions of~$V_1$ and~$V_2$ to contain at least one cluster in common, so we must declare the exploration as a failure. With this in mind, we define the second successful termination scenario as
\[A_2 := A_0' \cap A_1^c \cap \{\sigma < \infty\} \cap \{Z^* \text{ does not cause $V_1$ to increase}\}.\]
In case this occurs, we do not explore the region~$B_{\ell - \tau - \eta}(x) \cap \mathcal{R}_2 \cap \mathcal{R}_1^c$ further than~$\partial B_{\ell - \tau - \sigma}(x)$. 

Now assume that~$A_0' \cap A_1^c \cap \{\sigma = \infty\}$ occurs, that is,~$A_0' \cap A_1^c$ occurs and
\[\mathcal{X} \cap (B_{\ell - \tau - \eta}(x) \backslash B_{\ell - \tau - r_0 + \eta}(x)) \cap \mathcal{R}_2 \cap \mathcal{R}_1^c = \varnothing.\] We then explore the annulus~$B_{\ell - \tau}(x) \backslash B_{\ell - \tau - \eta - 1}(x)$, meaning that we reveal~$\mathcal{Y}_{\tau + \eta + 1}$. We let~$V_1'$ and~$V_2'$ be the extensions of~$V_1$ and~$V_2$, respectively, after this revealment. In other words,~$V_1'$ is the set of vertices of~$\mathcal{Y}_{\tau + \eta + 1}$ belonging to the simplices of any cluster of~$\mathcal{Y}_{\tau + \eta + 1}$ that intersect~$V_1$, and similarly for~$V_2'$. Then define
\[A_3 := A_0' \cap A_1^c \cap \{\sigma = \infty\} \cap \{V_1' \neq V_2'\} \cap \{V_2' \neq V_2\}.\]

This completes the description of the exploration, and we now turn to the proofs of~\eqref{eq_mutual} and~\eqref{eq_scenarios}.

\begin{proof}[Proof of~\eqref{eq_mutual}]
It is clear that~$\mathrm{Piv}^{1}_{n,\ell}(x)\cup \mathrm{Piv}^{2}_{n,\ell}(x) \subset A_0'$. 
We now argue that
\[(\mathrm{Piv}^{1}_{n,\ell}(x)\cup \mathrm{Piv}^{2}_{n,\ell}(x)) \cap A_1^c \cap A_2^c \subset A_3.\]
	Assume that the event on the left-hand side occurs. Then,~$A_0' \cap A_1^c$ also occurs, so we have~$V_2 \subset (B_{\ell - \tau + \eta}(x))^c$, which implies that
\begin{equation} \label{eq_implication_r2}\mathcal{R}_2 = \cup_{y \in V_2} B_{r_0}(y) \subset (B_{\ell - \tau - r_0 + \eta}(x))^c.\end{equation}
	Moreover,~$A_0' \cap A_1^c \cap \{\sigma = \infty\}$ also occurs, so 
\begin{equation*}
	\mathcal{X} \cap \mathcal{R}_2 \cap \mathcal{R}_1^c \cap B_{\ell - \tau - \eta}(x) \cap (B_{\ell - \tau - r_0 + \eta}(x))^c = \varnothing.
\end{equation*}
	Together with~\eqref{eq_implication_r2}, this gives
\begin{equation} \label{eq_implication_r22}
	\mathcal{X} \cap \mathcal{R}_2 \cap \mathcal{R}_1^c \cap B_{\ell - \tau - \eta}(x)  = \varnothing.
\end{equation}
Next, the occurrence of either of the pivotality events implies that there exists in~$\mathcal{X}$ a simplex that extends~$V_2$ to the interior of~$B_{\ell - \tau}(x)$ (while being distant from~$V_1$). More precisely, there is a simplex~$S \subset \mathcal{X}$ such that
\begin{equation*} 
S \cap B_{\ell - \tau}(x) \neq \varnothing,\qquad S \cap \mathcal{R}_2 \neq \varnothing,\qquad S \cap \mathcal{R}_1 = \varnothing.
\end{equation*}
In particular,~$S$ has a vertex~$w \in \mathcal{X} \cap \mathcal{R}_2 \cap \mathcal{R}_1^c$. Now,~\eqref{eq_implication_r22} implies that~$w \in (B_{\ell - \tau - \eta}(x))^c$. We then have~$S \subset B_1(w) \subset (B_{\ell - \tau - \eta - 1}(x))^c$. This implies that~$A_3$ occurs, so the proof of~\eqref{eq_mutual} is complete.
\end{proof}

It remains to prove~\eqref{eq_scenarios} in each of the three success scenarios. In each case, we will condition on the region that has already been explored, and describe a certain event on the unexplored region which occurs with probability higher than some~$\varepsilon_1 > 0$, and guarantees that~$\mathrm{Piv}^2_{n,\ell}(x)$ occurs. The construction is very similar to that of~\citep{strictIneq}. 
\begin{proof}[Proof of~\eqref{eq_scenarios}, case~$j = 1$]
	Figure~\ref{fig_eventA1} should aid in the understanding of the construction that follows.
	Recall that~$Y$ is the vertex of~$V_1$ in~$\partial B_{\ell - \tau}(x)$. When~$A_1$ occurs, there is also a vertex~$Z \in V_2 \cap B_{\ell - \tau + \eta}(x)$. We claim that, if~$\ell$ is large enough (and hence~$\ell - \tau$ is also large) and~$\eta$ is small enough, we can find points~$Y',Z' \in \R^2$ and~$\delta \in (0,1)$ such that
	\[ B_\delta(Y') \subset B_{r_0}(Y) \cap B_{\ell - \tau}(x),\quad B_\delta(Z') \subset B_{r_0}(Z) \cap B_{\ell - \tau}(x), \quad \vert Y' - Z'\vert  > r_0 + 2\delta\]
	and additionally,
	\begin{equation}\label{eq_want_ineq1} \mathrm{dist}(Y',(B_{\ell - \tau}(x))^c) > 1 + \delta, \quad \mathrm{dist}(Y',V_2) \ge \mathrm{dist}(Y', (B_{\ell - \tau}(x))^c \cap (B_{r_0}(Y))^c) > r_0 + \delta\end{equation}
	and
	\begin{equation}\label{eq_want_ineq2} \mathrm{dist}(Z',(B_{\ell - \tau}(x))^c) > 1 + \delta, \quad \mathrm{dist}(Z',V_1) \ge \mathrm{dist}(Z', (B_{\ell - \tau}(x))^c \cap (B_{r_0}(Z))^c) > r_0 + \delta.\end{equation}
	This claim, as well as similar claims that will be made in the proofs of the other cases, follows from elementary considerations of plane geometry, so we omit the details. In order to prove it, it is helpful to first pretend that~$\eta$ could be~$0$ (so that~$Z$ would belong to~$\partial B_{\ell - \tau}(x))$ and that~$\ell$ is so large that~$\partial B_{\ell - \tau}(x)$ can be replaced by a line, and then to obtain the desired conclusion in our actual situation by using continuity considerations.

	With~\eqref{eq_want_ineq1} and~\eqref{eq_want_ineq2} at hand, the following holds. Assume that, apart from the vertices already found in the exploration, we force~$\mathcal{X}$ to include~$d$ vertices inside~$B_\delta(Y')$ (and nothing more, for now). Then, the simplex formed by these vertices will be within distance~$r_0$ to~$Y$, so it will join~$V_1$, but it will stay at distance larger than~$r_0 + \delta$ from~$V_2$ and~$Z'$. Also, no additional simplex will be formed between the new vertices and previously existing ones (the distance between new and old vertices is more than~$1+\delta$). Similar considerations hold if instead, we force~$\mathcal{X}$ to include~$d$ vertices inside~$B_\delta(Z')$.

	Next, we can also obtain sequences of points~$y_0 ,y_1,\ldots,y_M$ and~$z_0, z_1,\ldots,z_N$ inside~$B_{\ell - \tau}(x)$ such that
	\begin{gather*}
		y_0 = Y',\quad z_0 = Z',\quad y_1,\ldots,y_M,z_1,\ldots,z_N \in B_{\ell - \tau - r_0}(x), \quad y_M, z_N \in B_{r_0}(x),\\[.2cm]
		\vert y_{a+1}-y_a \vert \le r_0 - 2\delta,\; a=0,\ldots, M-1,\\[.2cm] \vert z_{b+1}-z_b \vert \le r_0 - 2\delta,\; b=0,\ldots, N-1,\\[.2cm]
	\min\{\vert y_a - z_b \vert: 0 \le a \le M,\; 0 \le b \le N\} > r_0 +2\delta.
	\end{gather*}
\begin{figure}[htb]
\begin{center}
\setlength\fboxsep{0pt}
\setlength\fboxrule{0.5pt}
\fbox{\includegraphics[width = 0.7\textwidth]{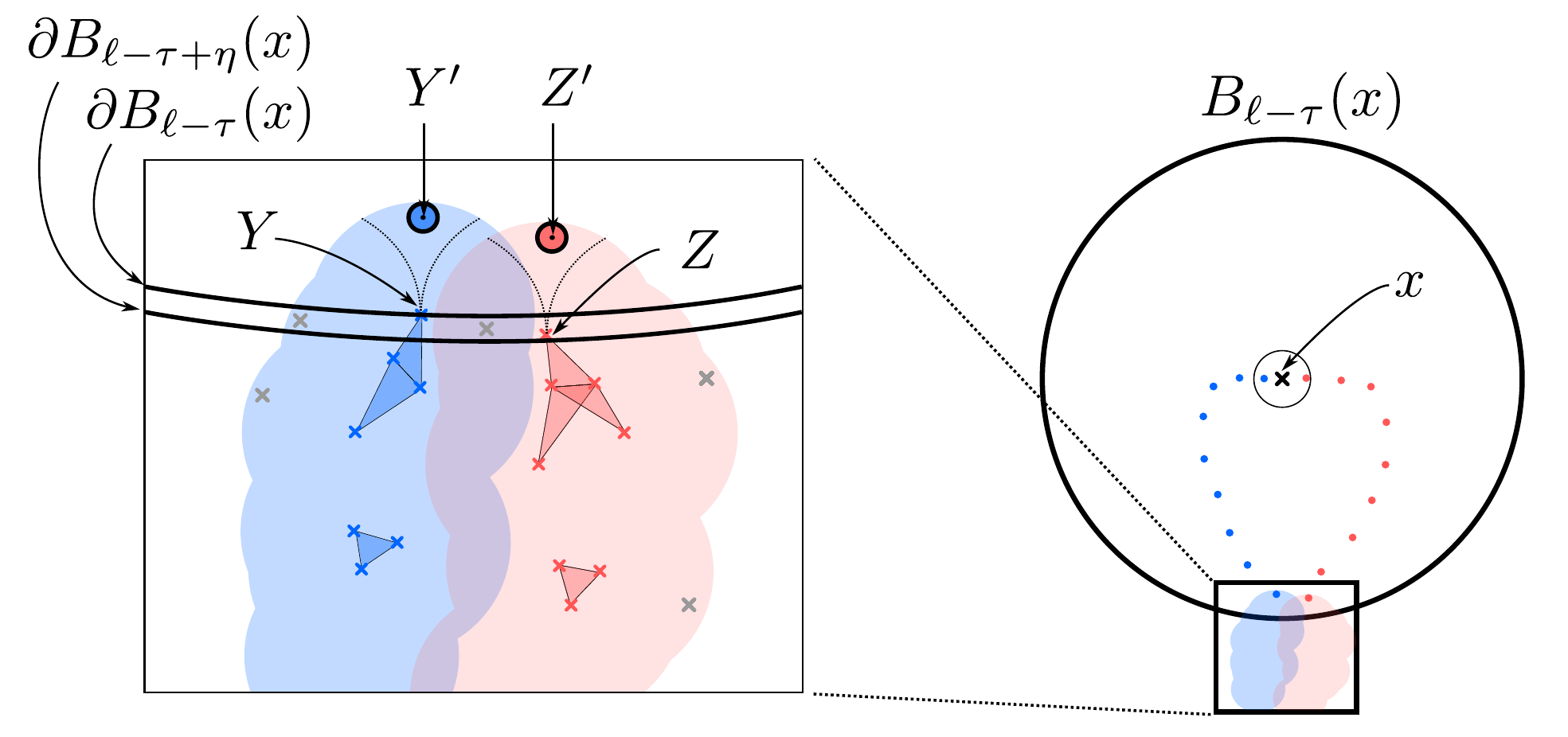}}
\end{center}
\caption{Summary of the construction of an event that guarantees pivotality in case~$A_1$ occurs.}
\label{fig_eventA1}
\end{figure}

	We now let~$\tilde{A}_1$ be the event that~$A_1$ occurs, that~$X$ contains exactly~$d$ vertices in each of the balls~$B_\delta(y_0),\ldots, B_{\delta}(y_M), B_\delta(z_0),\ldots, B_\delta(z_N)$, and no other vertex in~$B_{\ell - \tau}(x)$, and finally that~$\Xi_1 \equiv 1$ on~$(X \cap B_{\ell - \tau}(x)) \cup \{x\}$. We then have~$\mathrm{Piv}^2_{n,\ell}(x) \supset \tilde{A}_1$
	and~$\mathbb{P}(\tilde{A}_1 \mid \mathcal{G}) \cdot \mathds{1}_{A_1} \ge \varepsilon_1 \cdot \mathds{1}_{A_1}$ for some~$\varepsilon_1 > 0$ depending only on~$\ell$,~$p$ and~$q$.
\end{proof}

\begin{figure}[htb]
\begin{center}
\setlength\fboxsep{0pt}
\setlength\fboxrule{0.5pt}
\fbox{\includegraphics[width = 0.7\textwidth]{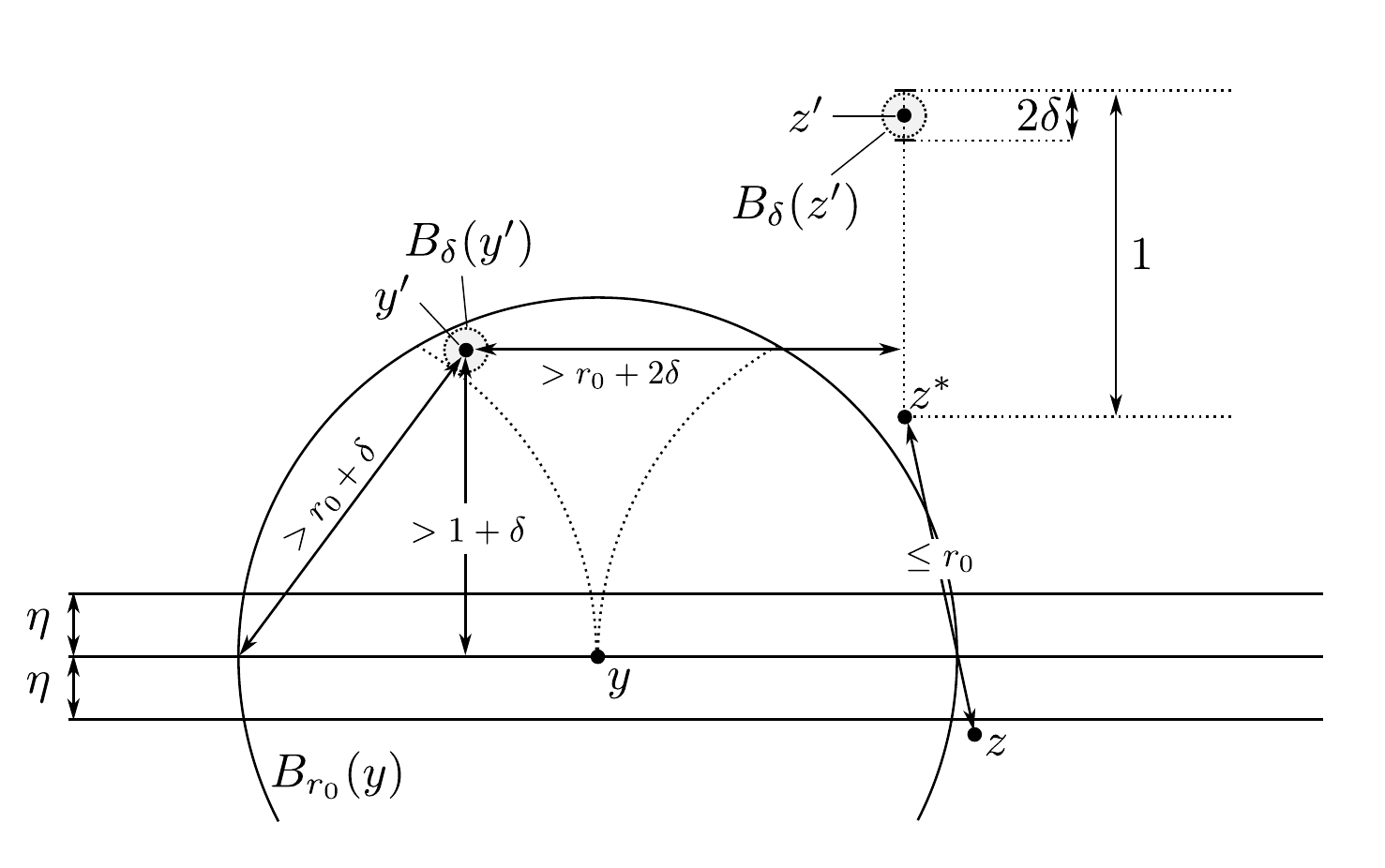}}
\end{center}
\caption{Configuration of points for the treatment of the event~$A_2$}
\label{fig_eventA2}
\end{figure}

Before turning to the proof of~\eqref{eq_scenarios} for~$j=2$, we describe a configuration of points in~$\R^2$ which will later be mapped to points near~$\partial B_{\ell - \tau}(x)$ in the event~$A_2$. See Figure~\ref{fig_eventA2}. Let~$z = (z_1,z_2)$,~$z^* = (z^*_1,z^*_2)$,~$y=(y_1,y_2) \in \R^2$ be three points satisfying
	\[y_2 = 0,\quad  z_2 < -\eta, \quad z_2^* > \eta,\quad \vert z - y\vert  > r_0,\quad \vert z'-y\vert  > r_0,\quad\vert z-z^*\vert  < r_0.\]
	Then, recalling that~$r_0 > 2$, it is easy to check that there exists~$y' = (y'_1,y'_2)$ and~$\delta > 0$ such that
	\begin{gather*} B_\delta(y') \subset B_{r_0}(y),\quad \mathrm{dist}(y',(B_{r_0}(y))^c \cap (\R \times (-\infty,0])) > r_0 + \delta,\\ 
	y'_2 > 1+ \delta, \quad \vert  y'_1- z^*_1\vert  > r_0+2\delta.\end{gather*}
	We assume that~$\delta < \eta/2$.
	Now let~$z'= (z'_1,z'_2) := (z^*_1,z^*_2 + 1 - \delta)$. We have
	\[\max_{a \in B_\delta(z')} \vert a-z^*\vert  < 1,\quad \vert y' - z'\vert  > r_0 + 2 \delta,\quad z_2' > 1+\delta,\]
	the last inequality following from~$z^*_2 > \eta > 2\delta$.

\begin{proof}[Proof of~\eqref{eq_scenarios}, case~$j = 2$]
	The main difference from the case~$j=1$ is the way in which the points~$Y'$ and~$Z'$ are chosen. Let us describe it. In the event~$A_2$, we succeed in finding a vertex~$Z^* \in \mathcal{X} \cap \mathcal{R}_2 \cap \mathcal{R}_1^c \cap B_{\ell - \tau - \eta}(x)$. 
	Since~$Z^* \in \mathcal{R}_2$, there exists some~$Z \in V_2$ such that~$\vert Z-Z^*\vert  \le r_0$, and since~$A_1$ does not occur, we have that~$Z \in (B_{\ell - \tau + \eta}(x))^c$. We now make the comparison with the configuration in the plane described above, with~$B_{\ell - \tau}(x)$ corresponding to the upper half plane, and~$Y,Z,Z^*,Y',Z'$ corresponding to~$y,z,z^*,y',z'$, respectively ($Z'$ is taken at distance~$1 - \delta$ from~$Z^*$ in the segment between~$Z^*$ and~$x$). If~$\ell$ is large enough, the choice of~$\delta$ is still possible and we have
	\begin{gather*}\vert Y' - Z^*\vert  > r_0 + 2\delta,\quad \vert Y' - Z'\vert  > r_0 + 2\delta,\\[.2cm]
		\mathrm{dist}(Y',(B_{\ell - \tau}(x))^c) > 1 + \delta, \\[.2cm] \mathrm{dist}(Y',V_2) \ge \mathrm{dist}(Y', (B_{\ell - \tau}(x))^c \cap (B_{r_0}(Y))^c) > r_0 + \delta,\\[.2cm]
		\mathrm{dist}(Z',(B_{\ell - \tau}(x))^c) > 1 + \delta,\quad \max_{a \in B_{\delta}(Z')}\vert a-Z^*\vert  < 1.
	\end{gather*}
	Next, note that the distance to~$V_1$ increases strictly as we move from~$Z^*$ to~$Z'$, so (decreasing~$\delta$ if necessary) we also have
	\[\mathrm{dist}(Z',V_1) > \mathrm{dist}(Z^*,V_1) + \delta > r_0 + \delta;\]
	the second inequality holds since~$Z^* \in \mathcal{R}_1^c$.

	We now complete the argument similarly to the case~$j = 1$. We choose vertices~$y_0 = Y', y_1,\ldots, y_M$,~$z_0 = Z', z_1,\ldots,z_N$ as before, and consider the event that, apart from the vertices already seen in the exploration,~${X}$ contains~$d$ vertices inside the balls of radius~$\delta$ centered at each of these vertices, and nothing more (in particular, the annulus~$B_{\ell - \tau}(x) \backslash B_{\ell - \tau - \eta}(x)$ stays empty) and that~$\Xi^1 \equiv 1$ on~$(X \cap B_{\ell - \tau}(x)) \cup \{x\}$. 
\end{proof}

\begin{figure}[htb]
\begin{center}
\setlength\fboxsep{0pt}
\setlength\fboxrule{0.5pt}
\fbox{\includegraphics[width = 0.7\textwidth]{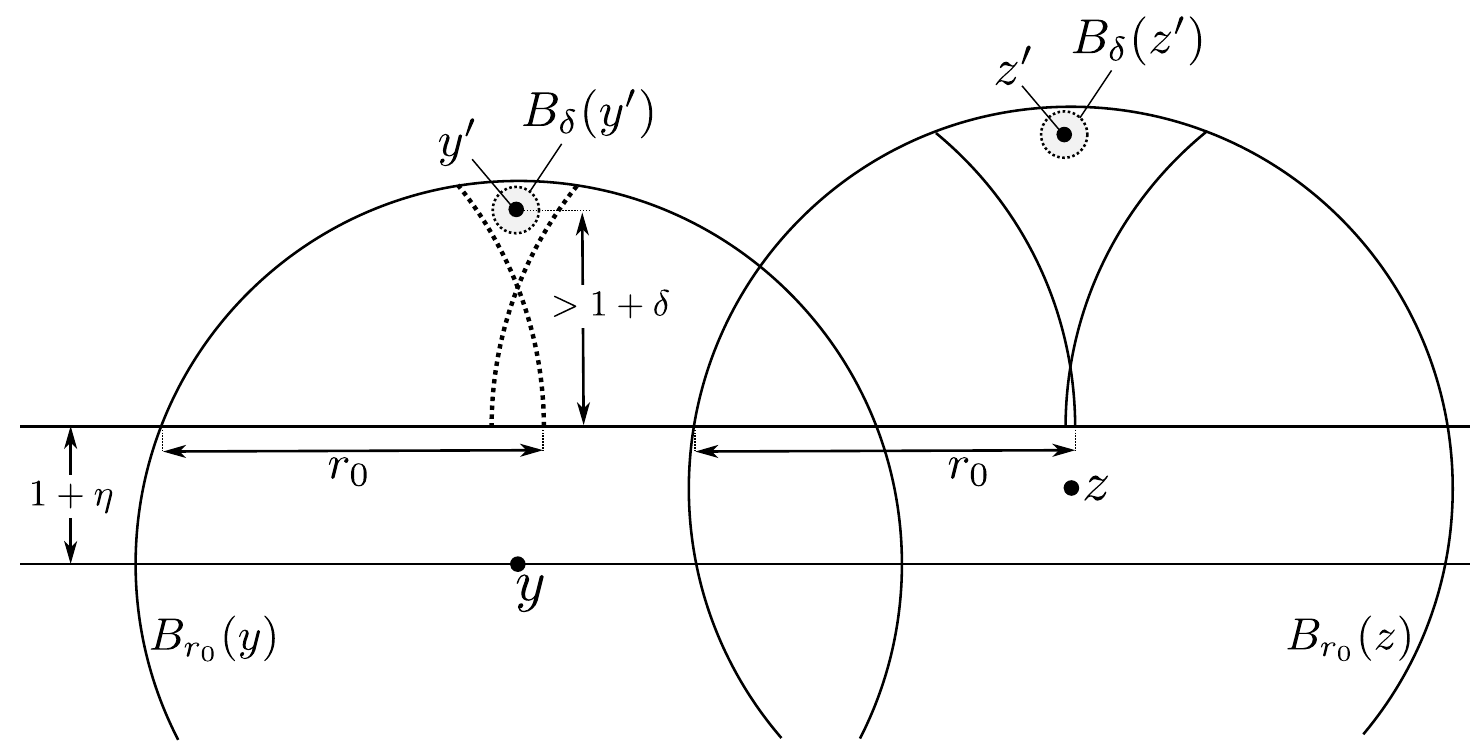}}
\end{center}
\caption{Configuration of points for the treatment of the event~$A_3$}
\label{fig_eventA3}
\end{figure}

The proof of~\eqref{eq_scenarios} for the case~$j = 3$ will again be preceded by a description of a configuration of points in the plane. See Figure~\ref{fig_eventA3}. We let~$y = (y_1,y_2),z = (z_1,z_2) \in \R^2$ be such that
\[y_2 = 0,\quad z_2 \in (0, 1+\eta),\quad \vert z-y\vert \ge  r_0.\]
Recalling that~$r_0 > 2$, if~$\eta$ and~$\delta$ are small enough we can then find~$y' = (y'_1,y'_2),z'=(z'_1,z'_2)$ such that~$\vert y'-z'\vert  > r_0 + 2\delta$ and
\begin{gather*} B_\delta(y') \subset B_{r_0}(y),\quad  y'_2 > 2+\eta+\delta,\\ \mathrm{dist}(y', (B_{r_0}(y))^c \cap (\R \times (-\infty, 1+\eta])) > r_0 + \delta,\\[.2cm]
B_\delta(z') \subset B_{r_0}(z),\quad  z'_2 > 2+\eta+\delta,\\ \mathrm{dist}(z', (B_{r_0}(z))^c \cap (\R \times (-\infty, 1+\eta])) > r_0 + \delta
\end{gather*}

\begin{proof}[Proof of~\eqref{eq_scenarios}, case~$j = 3$]
	In the event~$A_3$, there is a vertex~$Z$ of~$V_2'$ in $B_{\ell-\tau}(x) \backslash B_{\ell-\tau-\eta-1}(x)$. We now make the correspondence with the configuration of points just described, as before, so that we find points~$Y'$ and~$Z'$ matching the requirements satisfied by~$y'$ and~$z'$. Note that this implies that vertices of~$\mathcal{X}$ introduced inside~$B_\delta(Y')$ will not form a simplex with vertices outside~$B_{\ell - \tau - \eta - 1}(x)$, and a simplex entirely contained in~$B_\delta(Y')$ will be at distance larger than~$r_0$ from~$V_2'$, and also from~$B_\delta(Z')$. Similar considerations apply to~$Z'$. The proof is now completed as in the other cases.
\end{proof}

\subsection*{Acknowledgement.}
 The authors are indebted to D.~Yogeshwaran for making us aware of the existing works on phase transitions and sharp thresholds in the \v Cech complex, and also of the relation with clique percolation.

\bibliography{lit.bib}

\end{document}